\documentclass[11pt]{article}
\usepackage{amsmath,amssymb,mathrsfs}
\usepackage{amsfonts}
\usepackage{color}
\newtheorem{theo}{Theorem}[section]
\newtheorem{lem}[theo]{Lemma}
\newtheorem{cor}[theo]{Corollary}
\newtheorem{rem}[theo]{Remark}
\newtheorem{prop}[theo]{Proposition}

\newcommand{\mysection}[1]{\section{#1} \setcounter{equation}{0}}
\newcommand{\proof}{{\sc Proof.} \quad}
\newcommand{\proofc}{{\sc Proof} \ }
\newcommand{\be}{\begin{equation} \label}
\newcommand{\ee}{\end{equation}}
\newcommand{\bea}{\begin{eqnarray}\label}
\newcommand{\eea}{\end{eqnarray}}
\newcommand{\bas}{\begin{eqnarray*}}
\newcommand{\eas}{\end{eqnarray*}}
\newcommand{\bit}{\begin{itemize}}
\newcommand{\eit}{\end{itemize}}
\newcommand{\qed}{\hfill$\Box$ \vskip.2cm}
\newcommand{\nn}{\nonumber}
\newcommand{\R}{\mathbb{R}}
\newcommand{\N}{\mathbb{N}}

\newcommand{\eps}{\varepsilon}

\newcommand{\abs}{\\[5pt]}

\newcommand{\uv}{\underline{v}}
\newcommand{\ou}{\overline{u}}

\newcommand{\oy}{\overline{y}}

\newcommand{\hc}{\widehat{C}}
\newcommand{\hhc}{\widetilde{C}}

\newcommand{\f}{f_A^{(\alpha,\beta)}}
\newcommand{\xis}{\xi_\star}
\newcommand{\ua}{u_A^{(\alpha)}}
\newcommand{\tu}{\widetilde{u}}
\begin{document}
\title{Slow growth of solutions of super-fast diffusion equations\\
with unbounded initial data}
\author{
Marek Fila\footnote{fila@fmph.uniba.sk}\\
{\small Department of Applied Mathematics and Statistics, Comenius University,}\\
{\small 84248 Bratislava, Slovakia}
\and
Michael Winkler\footnote{michael.winkler@math.uni-paderborn.de}\\
{\small Institut f\"ur Mathematik, Universit\"at Paderborn,}\\
{\small 33098 Paderborn, Germany} }
\date{}
\maketitle
\begin{abstract}
\noindent 
We study positive solutions of the super-fast diffusion equation in the
whole space with initial data which are unbounded as $|x|\to\infty$. We find
an explicit dependence of the slow temporal growth rate of solutions on the
initial spatial growth rate. A new class of self-similar solutions plays a
significant role in our analysis.
\abs
\end{abstract}
\newpage
\mysection{Introduction}
Investigating mechanisms of mass flux plays an important role in the literature on nonlinear diffusion processes.
In the specific context of the equation		
\be{NLD}
	v_t=\nabla \cdot (v^{m-1}\nabla v)
\ee
when posed in the entire space $\R^n$, by a large number of results quite a comprehensive understanding 
has been achieved with regard to phenomena related to mass {\em outflux}.
For instance, in the case $m>1$ when (\ref{NLD}) becomes the porous medium equation, and also within the range
$(n-2)_+/n<m<1$ of fast diffusion,
the large time behaviour of nonnegative 
solutions which decay sufficiently fast at spatial infinity, in the sense of having finite mass $\int_{\R^n} v$,
is essentially determined by a particular element of a family of explicit 
self-similar solutions, the so-called Barenblatt solutions (\cite{friedman_kamin}, \cite{vazquez_JEE2003},
\cite{vazquez_book}). For solutions of (\ref{NLD}) with finite initial mass, the mass is conserved if
$m\ge (n-2)/n$ when $n>2$ and $m>0$ when $n=1,2$ (\cite{vazquez_book}). If $n=1$ and $-1<m\le 0$ then there
is nonuniqueness and conservation of mass holds for the maximal solution
(\cite{ERV}, \cite{vazquez_book}). In the remaining cases $n=1,2$ and $m\le
(n-2)_+$ or $n>2$ and 
$m< (n-2)/n$, the mass is not conserved
(\cite{vazquez_book}).
In the borderline case $m=(n-2)/n$,
$n>2$,
the large time behaviour of solutions with finite mass is more complicated than for
$m>(n-2)/n$ because the solution does not evolve toward a single
self-similar solution. The behaviour is different in an inner region and in
an outer one (\cite{GPV}, \cite{King}).\abs
In the case $m<(n-2)/n$, mass flux toward spatial infinity occurs in an
effective manner: If $0<m<(n-2)/n$ 
then any positive solution emanating from initial data in $L^{n(1-m)/2}(\R^n)$ becomes extinct in finite time,
and in the super-fast diffusion range $m<0$ even instantaneous extinction occurs for such initial data in the sense
that then no local-in-time solution exists (\cite{DP97}, cf.~also \cite{DP99}
for the case $m=0$).
Beyond this, the literature has provided more detailed information on how the asymptotic
behaviour near extinction 
depends on the initial spatial decay, again indicating an important role
of self-similar solutions (see e.g.~\cite{BBDGV}, \cite{BDGV}, \cite{BGV},
\cite{CV}, 
\cite{fvwy_ARMA}, \cite{fkw_JLMS}, \cite{GP}, \cite{King}, \cite{vazquez_book}).\abs
In contrast to this, only little seems known about processes of mass {\em influx} from
infinity,
except for few results on essentially one-dimensional wave-like transport mechanisms (\cite{aronson_lect}; 
cf.~also \cite{win_MATANN} and \cite{win_JDDE} for two recent examples involving non-constant wave speeds).
In the present work we establish some results in this direction by exploring in detail how `mass'
initially concentrated at spatial infinity spreads over the entire space for
the super-fast diffusion equation (\ref{NLD}) with $m<0$. 
More precisely, we shall be concerned with the Cauchy problem
\be{0fast}
	\left\{ \begin{array}{ll}
	v_t=\nabla \cdot (v^{m-1} \nabla v), \qquad &x\in\R^n, \ t>0, \\[1mm]
	v(x,0)=v_0(x), \qquad &x\in\R^n,
 	\end{array} \right.
\ee
in the strongly degenerate regime $m<0$, assuming the initial data $v_0\in C^0(\R^n)$ are positive,
and such that 
\be{growth}
	v_0(x)\to + \infty
	\qquad \mbox{as } |x|\to\infty
\ee
in an appropriate sense.
Then solutions exist globally (\cite{DP97})  
and it is natural to expect that they tend to $+\infty$ everywhere in $\R^n$ in the large time limit. 
Our main objective now consists in investigating quantitatively the dependence of such growth phenomena on the
particular asymptotics of the initial data. Let us mention
here that for $m<0$ the diffusion is very fast where $v$ is small but very
slow where $v$ is large. Therefore, it is natural to expect a slow growth
process.\abs
{\bf Main results: Decay estimates in a degenerate parabolic equation.}\quad
In order to transfer the above situation to a convenient framework involving bounded functions decaying at spatial
infinity, we substitute 
\be{subst}
	u(x,t):=v^m(x,t), \qquad x\in\R^n, \ t>0,
\ee
and then obtain formal equivalence of (\ref{0fast}) to the Cauchy problem
\be{0}
	\left\{ \begin{array}{ll}
	u_t=u^p \Delta u, \qquad &x\in\R^n, \ t>0, \\[1mm]
	u(x,0)=u_0(x), \qquad &x\in\R^n,
 	\end{array} \right.
\ee
where $p:=(m-1)/m=(|m|+1)/|m|>1$ and $u_0:=v_0^m$ is positive and continuous in $\R^n$.
Actually, most parts of our analysis also apply to the case $p=1$ which does
not stem from super-fast diffusion.\abs
In view of known results on nonuniqueness of classical solutions to (\ref{0fast}) 
(\cite{RV}), even in the framework of smooth positive solutions we cannot expect solutions of (\ref{0}) 
to be uniquely determined. 
As a preliminary to our subsequent analysis, we shall therefore first make sure that after all, 
(\ref{0}) possesses a {\em minimal} classical solution for any positive continuous and bounded initial data.
\begin{prop}\label{prop201}
  Let $p\ge 1$, and suppose that $u_0\in C^0(\R^n) \cap L^\infty(\R^n)$ is positive.
  Then problem (\ref{0}) possesses at least one global classical solution $u\in C^0(\R^n\times [0,\infty))
  \cap C^{2,1}(\R^n\times (0,\infty))$ which 
  satisfies 
  \be{201.1}
	0<u(x,t)\le \|u_0\|_{L^\infty(\R^n)}	
	\qquad \mbox{ for all $x\in\R^n$ and } t\ge 0.
  \ee
  Moreover, this solution is minimal in the sense that whenever $T\in (0,\infty]$ and
  $\tu \in C^0(\R^n\times [0,T)) \cap C^{2,1}(\R^n\times (0,T))$ are such that $\tu$ is positive and
  solves (\ref{0}) classically
  in $\R^n\times (0,T)$, we necessarily have $u\le \tu$ in $\R^n\times (0,T)$.
\end{prop}
Now if $u_0$ belongs to $L^{q_0}(\R^n)$ for some positive $q_0$ not necessarily exceeding the value $1$,
we can establish the following implications for the temporal decay of
$\|u(\cdot,t)\|_{L^q(\R^n)}$
for $q\in (q_0,\infty]$.
\begin{theo}\label{theo200}
  Let $p\ge 1$, and suppose that $u_0\in C^0(\R^n) \cap L^\infty(\R^n)$ is positive and such that
  $u_0\in L^{q_0}(\R^n)$ for some $q_0>0$.
  Then for any $q>q_0$ one can find $C=C(q)>0$ with the property that the solution $u$ of (\ref{0}) 
  from Proposition~\ref{prop201} satisfies
  \be{200.1}
	\|u(\cdot,t)\|_{L^q(\R^n)} \le C  t^{-(1-\frac{q_0}{q})/(p+\frac{2q_0}{n})}
	\qquad \mbox{for all } t>0.
  \ee
  Moreover, for any $\delta>0$ there exists $\hc(\delta)>0$ such that
  \be{200.2}
	\|u(\cdot,t)\|_{L^\infty(\R^n)} \le \hc(\delta) t^{-n/(np+2q_0)+\delta}
	\qquad \mbox{for all } t>0.
  \ee
  In particular, if $u_0\in \bigcap_{q_0>0} L^{q_0}(\R^n)$ then for any $\delta>0$ one can
  find
  $\hhc(\delta)>0$ such that
  \be{200.3}
	\|u(\cdot,t)\|_{L^\infty(\R^n)} \le \hhc(\delta) t^{-\frac{1}{p}+\delta}
	\qquad \mbox{for all } t>0.
  \ee
\end{theo}
A natural next question appears to be how far the above one-sided estimates are optimal.
Surprisingly, the above result on decay in $L^\infty(\R^n)$ for initial data with fast decay cannot
be substantially improved in the sense that not even choosing $\delta=0$ is possible in (\ref{200.3}):
\begin{prop}\label{prop103}
  Let $p\ge 1$. Then for every positive $u_0\in C^0(\R^n)$, any global positive classical $u$ of (\ref{0}) 
  has the property that for any $R>0$ we have
  \bas
	\inf_{|x|<R} \Big\{ t^\frac{1}{p} u(x,t) \Big\} \to +\infty
	\qquad \mbox{as } t\to\infty.
  \eas
\end{prop}
But, indeed, also (\ref{200.1}) and (\ref{200.2}) are essentially sharp:
\begin{theo}\label{theo100}
  Let $p\ge 1$ and $q_0>0$. Then for every $q \in (q_0,\infty]$ there exists a positive function
  $u_0\in C^0(\R^n) \cap L^\infty(\R^n)$ such that $u_0\in L^{q_0}(\R^n)$ and such that for any
  $\delta>0$ one can find $C(\delta)>0$ with the property that the solution $u$ of (\ref{0})
  from Proposition~\ref{prop201} satisfies
  \be{100.1}
	\|u(\cdot,t)\|_{L^q(\R^n)} \ge C(\delta)  t^{-(1-\frac{q_0}{q})/(p+\frac{2q_0}{n})-\delta}
	\qquad \mbox{for all } t>1.
  \ee
\end{theo}
In cases when the initial data satisfy pointwise decay estimates of algebraic type, we can even achieve more
precise information on the respective large time behaviour. 
Fundamental for our analysis in this direction will be the observation that at least for $p>1$,
(\ref{0}) possesses a two-parameter family of self-similar solutions with suitable spatial decay.
As these solutions apparently have not yet been detected anywhere in the literature, let us describe them
in the following separate statement.
\begin{theo}\label{theo1000}
  Let $p>1$, $\alpha\in (0,\frac{1}{p})$ and $A>0$. 
  Then the equation $u_t=u^p\Delta u$
  possesses a radially symmetric positive classical self-similar solution $\ua$ which can be written in the form
  \be{ua}
	\ua(x,t)=t^{-\alpha} f \Big( t^{-\beta} |x| \Big), \qquad x\in\R^n, \ t>0,
  \ee
  where $\beta:=\frac{1-p\alpha}{2}>0$, and where $f:[0,\infty)\to (0,\infty)$ is
  the solution of the initial value problem
\be{0f}
        \left\{ \begin{array}{l}
        f^p(\xi)  \Big( f''(\xi)+\frac{n-1}{\xi}f'(\xi) \Big) + \beta \xi f'(\xi)+ \alpha f(\xi)=0, \qquad \xi>0, \\[1mm]
        f(0)=A, \quad f'(0)=0.
        \end{array} \right.
\ee
Moreover, $f$ satisfies 
  \be{f_upper_lower}
	c (1+\xi)^{-\frac{\alpha}{\beta}} \le f(\xi) \le C (1+\xi)^{-\frac{\alpha}{\beta}}
	\qquad \mbox{for all } \xi\ge 0
  \ee
  with appropriate positive constants $c$ and $C$.
\end{theo}
By means of two arguments based on parabolic comparison, in the upper estimate 
involving suitable members of the above self-similar family and in the lower estimate relying on 
certain compactly supported separated solutions, it is possible to give quite a comprehensive description
of the temporal asymptotics in (\ref{0}) for algebraically decaying initial data:
\begin{theo}\label{theo2000}
  Let $p\ge 1$ and $u_0\in C^0(\R^n)\cap L^\infty(\R^n)$ be positive, and let $u$ denote the solution
  of (\ref{0}) constructed in Proposition~\ref{prop201}.\abs
  (i) \ If 
  \be{2000.3}
	u_0(x) \ge C_0 (1+|x|)^{-\gamma}
	\qquad \mbox{for all } x\in \R^n
  \ee
  with some $C_0>0$ and $\gamma>0$,
  then for any $q\in (\frac{n}{\gamma},\infty]$ there exists $C>0$ such that
  \be{2000.4}
	\|u(\cdot,t)\|_{L^q(\Omega)} \ge C t^{-(\gamma-\frac{n}{q})/(p\gamma+2)}
	\qquad \mbox{for all } t>1.
  \ee
  (ii) \ If $p>1$ and there exist $\gamma>0$ and $C_1>0$ such that
  \be{2000.1}
	u_0(x) \le C_1  (1+|x|)^{-\gamma}
	\qquad \mbox{for all } x\in\R^n,
  \ee
  then for any $q\in (\frac{n}{\gamma},\infty]$ 
  one can find $C>0$ with the property that
  \be{2000.2}
	\|u(\cdot,t)\|_{L^q(\R^n)} 
	\le C t^{-(\gamma-\frac{n}{q})/(p\gamma+2)}
	\qquad \mbox{for all } t\ge 0.
  \ee
\end{theo}
{\bf Main results: Growth estimates for (\ref{0fast})}. \quad
Among the obvious translations of the above results to the original problem, let us highlight some
implications for the respective maximal classical solutions $v:=u^{1-p}$, $p=(m-1)/m$, of (\ref{0fast})
obtained from Proposition~\ref{prop201}.
First, Theorem~\ref{theo200}, Proposition~\ref{prop103} and Theorem~\ref{theo100} immediately yield the following.
\begin{theo}\label{theo200f}
  Let $m<0$, and let $v_0\in C^0(\R^n)$ be positive and such that
  $v_0^m\in L^\infty(\R^n) \cap L^{q_0}(\R^n)$ for some $q_0>0$.
  Let $v$ be the maximal solution of (\ref{0fast}).
  Then for any $\delta>0$ there exists $C(\delta)>0$ such that
  \be{200.2f}
	\inf_{x\in\R^n} v(x,t) \ge C(\delta) t^{1/(1-m-\frac{2q_0 m}{n}) - \delta}
	\qquad \mbox{for all } t>0.
  \ee
  Moreover, if $u_0^m\in \bigcap_{q_0>0} L^{q_0}(\R^n)$ then for any $\delta>0$ one can find $\hhc(\delta)>0$ such that
  \be{200.3f}
	\inf_{x\in\R^n} v(x,t) \ge C(\delta) t^{1/(1-m-\frac{2q_0 m}{n}) - \delta}
	\qquad \mbox{for all } t>0.
  \ee
  On the other hand, 
  for every $q_0>0$ there exists a positive function
  $v_0\in C^0(\R^n)$ such that $v_0^m \in L^\infty(\R^n) \cap  L^{q_0}(\R^n)$ and such that for each
  $\delta>0$ one can find $\hc(\delta)>0$ with the property that 
  \be{100.1f}
	\inf_{x\in\R^n} v(x,t) \le \hc(\delta) t^{1/(1-m-\frac{2q_0 m}{n}) + \delta}
	\qquad \mbox{for all } t>0.
  \ee
  Moreover, if $v_0\in C^0(\R^n)$ is an arbitrary positive function such that
  $1/v_0\in L^\infty(\R^n)$,
  then
  \be{200.4f}
	t^{-\frac{1}{1-m}} v(\cdot,t) \to 0
	\qquad \mbox{in } L^\infty_{loc}(\R^n)
	\qquad \mbox{as } t\to\infty.
  \ee
\end{theo}
Next, assuming algebraic growth of the initial data, we may rephrase
Theorem~\ref{theo2000} as follows.
\begin{theo}\label{theo2000f}
  Let $m<0$ and $v_0\in C^0(\R^n)$ be positive.
  Let $v$ be the maximal solution of (\ref{0fast}).\abs
  (i) \ If 
  \be{2000.3f}
	v_0(x) \le C_0 (1+|x|)^{\theta}
	\qquad \mbox{for all } x\in \R^n
  \ee
  with some $C_0>0$ and $\theta>0$,
  then there exists $C>0$ such that
  \be{2000.4f}
	\inf_{x\in\R^n} v(x,t) \le C t^\frac{\theta}{(1-m)\theta+2}
	\qquad \mbox{for all } t>1.
  \ee
  (ii) \ If there exist $\theta>0$ and $C_1>0$ such that
  \be{2000.1f}
	v_0(x) \ge C_1 (1+|x|)^{\theta}
	\qquad \mbox{for all } x\in \R^n
  \ee
  then
  \be{2000.2f}
	\inf_{x\in\R^n} v(x,t) \ge C t^\frac{\theta}{(1-m)\theta+2}
	\qquad \mbox{for all } t>1.
  \ee
  with some $C>0$.
\end{theo}
\begin{rem}\label{remf}
Theorem~\ref{theo2000f} describes explicitly how slow the growth process is
and how it slows down as $m\to -\infty$. Namely, for $m<0$ and $\theta>0$ set
\[
\vartheta(\theta,m) :=\frac{\theta}{(1-m)\theta+2},
\]
which is the exponent from (\ref{2000.4f}), (\ref{2000.2f}). Then for every
fixed $\theta>0$ we have that $\vartheta(\theta,m)\to 0$ as $m\to -\infty$,
$\vartheta$ is increasing in both variables, and $0<\vartheta(\theta,m)<1$ 
for all $m<0$ and $\theta>0$.\abs
The mechanism of mass influx from infinity is completely different for the
linear heat equation. For example, for any positive even integer $k$, the
solution of the problem
\[
\left\{ \begin{array}{ll}
        u_t=u_{xx} \, , \qquad &x\in\R, \ t>0,
\\[1mm]
        u(x,0)=x^k, \qquad &x\in\R,
        \end{array} \right.\cdot
\]
is the heat polynomial
\[
H_k(x,t):=\sum_{i=0}^{k/2}\frac{k!}{i!(k-2i)!}x^{k-2i}t^i,
\]
and
\[
\inf_{x\in\R} H_k(x,t)=\frac{k!}{(k/2)!}t^{k/2},
\]
so the growth rate tends to infinity as $k\to\infty$ while
$\vartheta(\theta,m)<1/(1-m)$ for all $\theta >0$.
\end{rem}
The paper is organized as follows. We give a sketch of the proof of Proposition~\ref{prop201} in Section~2. The upper bounds from Theorem~\ref{theo200} are established in Section~3. We derive lower bounds used in the proofs of Proposition~\ref{prop103}, Theorem~\ref{theo100} and Theorem~\ref{theo2000}~(i) in Section~4. The proof of Theorem~\ref{theo100} is finished in Section~5. Self-similar solutions are studied in Section~6 and they are used there to prove Theorem~\ref{theo2000}~(ii).

\mysection{Global existence via approximation. Proof of Proposition~\ref{prop201}}
If $p>1$ then Proposition~1.1 follows from \cite{DP97}. Since we include
also the case $p=1$, we give a brief sketch of a proof which works for $p\ge 1$.\abs
In order to construct solutions to (\ref{0}), for $B_R:=\{|x|<R\}$, $R>0$ we consider the approximate problems
\be{0R}
	\left\{ \begin{array}{ll}
	u_{Rt}=u_R^p \Delta u_R, \qquad &x\in B_R, \ t>0, \\[1mm]
	u_R(x,t)=0, \qquad &x\in\partial B_R, \ t>0, \\[1mm]
	u_R(x,0)=u_{0R}(x), \qquad &x\in B_R,
 	\end{array} \right.
\ee
where $u_{0R} \in C^3(\bar B_R)$ satisfies $0<u_{0R}<u_0$ in $B_R$ and $u_{0R}=0$ on $\partial B_R$ as well as
\be{conv_R}
	u_{0R} \nearrow u_0 \quad \mbox{in $\R^n$ \qquad as } R\nearrow \infty.
\ee
Moreover, for $\eps\in (0,1)$ we consider
\be{0Reps}
	\left\{ \begin{array}{ll}
	u_{R\eps t}=u_{R \eps}^p \Delta u_{R\eps}, \qquad &x\in B_R, \ t>0, \\[1mm]
	u_{R\eps}(x,t)=\eps, \qquad &x\in\partial B_R, \ t>0, \\[1mm]
	u_{R\eps}(x,0)=u_{0R\eps}(x), \qquad &x\in B_R,
 	\end{array} \right.
\ee
where we have set
\bas
	u_{0R\eps}:=u_{0R}+\eps.
\eas
\begin{lem}\label{lem202}
  Let $p\ge 1$, and assume that $u_0\in C^0(\R^n) \cap L^\infty(\R^n)$ is positive.
  Then with $u_{0R}$ and $u_{0R\eps}$ as above, problem (\ref{0Reps}) possesses
for any $\eps\in (0,1)$ a global classical
  solution $u_{R\eps} \in C^0(\bar B_R \times [0,\infty)) \cap C^{2,1}(\bar B_R\times (0,\infty))$. 
  As $\eps\searrow 0$, we have that $u_{R\eps}\searrow u_R$ where  
  $u_R\in C^0(\bar B_R \times [0,\infty)) \cap C^{2,1}(B_R\times (0,\infty))$ 
  is a positive classical solution of (\ref{0R}). 
  Moreover, there exists a classical solution $u\in C^0(\R^n \times [0,\infty)) \cap C^{2,1}(\R^n\times (0,\infty))$
  of (\ref{0}) which is such that (\ref{201.1}) holds, and that 
  $u_R\nearrow u$ in $\R^n\times (0,\infty)$ as $R\nearrow \infty$.
\end{lem}
\proof
  Since all arguments are well-known, we may confine ourselves to sketching the main steps only
  and refer for details to \cite{win_cauchy}, for example.\\
  According to standard theory of quasilinear parabolic problems, each of these actually nondegenerate problems
  possesses a globally defined classical solution $u_{R\eps}$ which satisfies
  \be{ure}
	\eps \le u_{R\eps} \le \|u_{0R}\|_{L^\infty(B_R)} + \eps
	\qquad \mbox{in } B_R \times (0,\infty).
  \ee
  By parabolic comparison it then follows that as $\eps\searrow 0$ we have $u_{R\eps} \searrow u_R$
  in $B_R\times (0,\infty)$ with some limit function $u_R$. According to the interior positivity properties of $u_{0R}$, it can be seen by comparison that
  $\inf_{\eps\in (0,1)} \inf_{(x,t)\in K\times [0,T]} u_{R\eps}(x,t)>0$ for any compact $K\subset B_R$ and every $T>0$,
  which combined with parabolic Schauder estimates (\cite{LSU}) shows that the convergence $u_{R\eps} \to 0$ 
  actually takes place in $C^0_{loc}(\bar B_R \times [0,\infty)) \cap C^{2,1}_{loc}(B_R \times [0,\infty))$,
  and that $u_R$ is a classical solution of (\ref{0R}) which due to (\ref{ure}) satisfies
  \be{ur}
	0< u_R \le \|u_{0R}\|_{L^\infty(B_R)}
	\qquad \mbox{in } B_R \times (0,\infty).
  \ee
  Now on the basis of this and the monotone approximation property (\ref{conv_R}), one more comparison argument
  asserts that $u_R \nearrow u$ in $\R^n\times (0,\infty)$ as $R\nearrow\infty$, where $u$ is a limit function which
  according to (\ref{ur}) clearly satisfies (\ref{201.1}).
  Again by means of parabolic regularity theory, the two-sided estimate in (\ref{201.1}) guarantees that actually 
  $u_R\to u$ in $C^0_{loc}(\R^n\times [0,\infty)) \cap C^{2,1}_{loc}(\R^n \times (0,\infty))$, and that hence
 $u$ solves (\ref{0}) classically.
\qed
\proofc of Proposition~\ref{prop201}.\quad
  Taking $u:=\lim_{R\to\infty} u_R$ as provided by Lemma~\ref{lem202}, in view of the latter we only need to show
  the claimed minimality property of $u$.
  Given $T\in (0,\infty]$ and a positive classical solution $\tu$ of (\ref{0}) in $\R^n\times (0,T)$ from the indicated class,
  however, by comparison we see that for every $R>0$ we have
  $u_R < \tu$ in $B_R \times [0,T)$ and that hence $u\le \tu$ in
  $\R^n\times [0,T)$.
\qed
\mysection{Upper decay estimates}
\subsection{Upper bounds for $u$ in $L^q(\R^n)$ for $q<\infty$}
The following elementary inequality will be used in Lemma~\ref{lem1}.
\begin{lem}\label{lem2}
  Let $\beta>0$. Then there exists $\kappa(\beta)>0$ such that
  \be{2.1}
	(a-b)_+^\beta \ge \kappa(\beta)  a^\beta - b^\beta
	\qquad \mbox{for all $a\ge 0$ and } b\ge 0.
  \ee
\end{lem}
\proof
  We take $c_1(\beta)>0$ such that
  \bas
	(A+B)^\beta \le c_1(\beta)  (A^\beta + B^\beta)
	\qquad \mbox{for all $A\ge 0$ and } B\ge 0.
  \eas
  Then given $a\ge 0$ and $b\ge 0$, we apply this to $A:=(a-b)_+$ and $B:=b$ to obtain using the nonnegativity of $\beta$
  that
  \bas
	a^\beta \le \Big((a-b)_+ +b \Big)^\beta \le c_1(\beta)  \Big((a-b)_+^\beta + b^\beta \Big).
  \eas
  Therefore, (\ref{2.1}) holds if we let $\kappa(\beta):=1/c_1(\beta)$, for instance.
\qed
Now the first main step toward our upper estimates for solutions is accomplished by means of a 
standard testing procedure applied to the approximate problems (\ref{0Reps}), followed by taking limits
properly.
\begin{lem}\label{lem1}
  Let $p\ge 1$ and $u_0\in C^0(\R^n)$ be positive with
  \bas
	u_0\in L^{q_0}(\R^n)
	\qquad \mbox{for some } q_0>0.
  \eas
  Then for any $q>q_0$ there exists $C(q)>0$ such that for any $R>0$, the solution $u_R$ of (\ref{0R}) satisfies
  \be{1.1}
	\|u_R(\cdot,t)\|_{L^q(B_R)} \le C(q) \, t^{-(1-\frac{q_0}{q})/(p+\frac{2q_0}{n})}
	\qquad \mbox{for all } t>0.
  \ee
%
\end{lem}
\proof
  Observing that $u_{R\eps}$ is smooth in $\bar B_R \times [0,\infty)$, for arbitrary
  $q>0$ we may test (\ref{0Reps}) by $u_{R\eps}^{q-1}$ and integrate by parts. 
  Using that $u_{R\eps}\ge \eps$ and that $u_{R\eps}=\eps$ on $\partial B_R$, we thereby obtain
  \bea{1.2}
	\frac{1}{q} \frac{d}{dt} \int_{B_R} u_{R\eps}^q
	&=& \int_{B_R} u_{R\eps}^{p+q-1} \Delta u_{R\eps} 
	= - (p+q-1) \int_{B_R} u_{R\eps}^{p+q-2} |\nabla u_{R\eps}|^2 
	+ \int_{\partial B_R} u_{R\eps}^{p+q-1} \frac{\partial u_{R\eps}}{\partial\nu} \nn\\
	&\le& - (p+q-1) \int_{B_R} u_{R\eps}^{p+q-2} |\nabla u_{R\eps}|^2 
	= - \frac{4(p+q-1)}{(p+q)^2} \int_{B_R} \Big|\nabla u_{R\eps}^\frac{p+q}{2} \Big|^2
  \eea
for all $t>0$.
  Since $p+q-1$ is positive for any $q>0$ due to the fact that $p\ge 1$, this implies that, in particular,
  \be{1.22}
	\|u_{R\eps}(\cdot,t)\|_{L^{q_0}(B_R)} \le \|u_{0R\eps}\|_{L^{q_0}(B_R)}
	\qquad \mbox{for all } t>0.
  \ee
  In order to proceed, let us fix $q>q_0$ and note that then since $\frac{2q}{p+q} \le 2$, 
  an application of the Gagliardo-Nirenberg inequality provides $c_1>0$ such that
  \be{1.3}
	\|\varphi\|_{L^\frac{2q}{p+q}(\R^n)} \le c_1 \|\nabla \varphi\|_{L^2(\R^n)}^a 
	\|\varphi\|_{L^\frac{2q_0}{p+q}(\R^n)}^{1-a}
	\qquad \mbox{for all } \varphi \in W^{1,2}(\R^n),
  \ee
  where $a\in (0,1)$ is determined by the relation
  \bas
	-\frac{n(p+q)}{2q} = \Big(1-\frac{n}{2}\Big)a - \frac{n(p+q)}{2q_0}(1-a),
  \eas
  that is, where
  \be{1.33}
	a=\left(\frac{n(p+q)}{2}\left(\frac{1}{q_0}-\frac{1}{q}\right)\right)\left(1-\frac{n}{2}+\frac{n(p+q)}{2q_0}\right)^{-1}.
  \ee
  Now for fixed $t>0$, we apply this to
  \bas
	\varphi(x):=\left\{ \begin{array}{ll}
	u_{R\eps}^\frac{p+q}{2}(x,t) - \eps^\frac{p+q}{2}, \qquad & x\in B_R, \\[1mm]
	0, & x\in \R^n \setminus B_R,
	\end{array} \right.
  \eas
  which indeed defines a function $\varphi$ belonging to $W^{1,2}(\R^n)$, because 
  $u_{R\eps}(\cdot,t) \in C^1(\bar B_R)$ satisfies $u_{R\eps} \ge \eps$ in $B_R$ and $u_{R\eps}|_{\partial B_R}=\eps$.\\
  Accordingly, (\ref{1.3}) shows that
  \be{1.4}
	\Big\| u_{R\eps}^\frac{p+q}{2}(\cdot,t)-\eps^\frac{p+q}{2}\Big\|_{L^\frac{2q}{p+q}(B_R)}
	\le c_1 \Big\|\nabla u_{R\eps}^\frac{p+q}{2}(\cdot,t)\Big\|_{L^2(B_R)}^a 
	\Big\| u_{R\eps}^\frac{p+q}{2}(\cdot,t)-\eps^\frac{p+q}{2}\Big\|_{L^\frac{2q_0}{p+q}(B_R)}^{1-a}
	\qquad \mbox{for all } t>0,
  \ee
  where again using that $u_{R\eps}\ge \eps$ and recalling (\ref{1.22}) we can estimate
  \bea{1.5}
	\Big\| u_{R\eps}^\frac{p+q}{2}(\cdot,t)-\eps^\frac{p+q}{2}\Big\|_{L^\frac{2q_0}{p+q}(B_R)}^{1-a}
	&\le& \Big\| u_{R\eps}^\frac{p+q}{2}(\cdot,t) \Big\|_{L^\frac{2q_0}{p+q}(B_R)}^{1-a} 
	= \|u_{R\eps}(\cdot,t)\|_{L^{q_0}(B_R)}^{\frac{p+q}{2q_0}(1-a)} \nn\\
	&\le& \|u_{0R\eps}\|_{L^{q_0}(B_R)}^{\frac{p+q}{2q_0}(1-a)}
	\qquad \mbox{for all } t>0.
  \eea
  Therefore, (\ref{1.4}) entails that
  \be{1.6}
	\int_{B_R} \Big|\nabla u_{R\eps}^\frac{p+q}{2}(\cdot,t)\Big|^2
	\le \Big\{ c_1 \|u_{0R\eps}\|_{L^{q_0}(B_R)}^{\frac{p+q}{2q_0}(1-a)} \Big\}^{-\frac{2}{a}}
	 \Big\| u_{R\eps}^\frac{p+q}{2}(\cdot,t)-\eps^\frac{p+q}{2}\Big\|_{L^\frac{2q}{p+q}(B_R)}^\frac{2}{a}
	\qquad \mbox{for all } t>0.
  \ee
  Here we apply Lemma~\ref{lem2} to see that with $\kappa(\cdot)$ as in (\ref{2.1}) we have
  \bas
	\Big\| u_{R\eps}^\frac{p+q}{2}(\cdot,t)&-&\eps^\frac{p+q}{2}\Big\|_{L^\frac{2q}{p+q}(B_R)}^\frac{2}{a}
	= \Bigg\{ \int_{B_R} \Big|u_{R\eps}^\frac{p+q}{2}(\cdot,t)-\eps^\frac{p+q}{2}\Big|^\frac{2q}{p+q} 
		\Bigg\}^\frac{p+q}{qa} \\
	&\ge& \Bigg\{ \int_{B_R} \bigg( \kappa \Big(\frac{2q}{p+q}\Big) u_{R\eps}^q(\cdot,t) - \eps^q \bigg) 
		\Bigg\}^\frac{p+q}{qa} 
	= \Bigg\{ \kappa\Big(\frac{2q}{p+q}\Big)  \int_{B_R} u_{R\eps}^q(\cdot,t) - |B_R|  \eps^q 
		\Bigg\}^\frac{p+q}{qa} \\
	&\ge& \kappa\Big(\frac{p+q}{qa}\Big) \bigg( \kappa\Big(\frac{2q}{p+q}\Big) \bigg)^\frac{p+q}{qa}
	 \bigg( \int_{B_R} u_{R\eps}^q(\cdot,t)\bigg)^\frac{p+q}{qa}
	- \Big( |B_R|  \eps^q \Big)^\frac{p+q}{qa}
  \eas
  for all $t>0$. As a consequence, from (\ref{1.6}) we obtain $c_2>0$ and $c_3(R)>0$ such that writing 
  $A_{R\eps}:=\|u_{0R\eps}\|_{L^{q_0}(B_R)}^{-(p+q)(1-a)/(aq_0)}$ we have
  \bas
	\int_{B_R} \Big|\nabla u_{R\eps}^\frac{p+q}{2}(\cdot,t)\Big|^2
	&\ge& c_2 A_{R\eps}  \bigg( \int_{B_R} u_{R\eps}^q(\cdot,t)\bigg)^\frac{p+q}{qa}
	- c_3(R) A_{R\eps} \eps^\frac{p+q}{a}
	\qquad \mbox{for all } t>0,
  \eas
  so that (\ref{1.2}) yields
  \bas
	\frac{d}{dt} \int_{B_R} u_{R\eps}^q(\cdot,t)
	\le -c_4 A_{R\eps}  \bigg(\int_{B_R} u_{R\eps}^q(\cdot,t)\bigg)^\frac{p+q}{qa}
	+ c_5(R) A_{R\eps} \eps^\frac{p+q}{a}
	\qquad \mbox{for all } t>0,
  \eas
  where again since $p\ge 1$, both $c_4:=\frac{4q(p+q-1)}{(p+q)^2} c_2$ and $c_5(R):=\frac{4q(p+q-1)}{(p+q)^2} c_3(R)$
  are positive. 
  By an ODE comparison argument, this entails that
  \be{1.7}
	\int_{B_R} u_{R\eps}^q(\cdot,t) \le \oy_{R\eps}(t)
	\qquad \mbox{for all } t>0,
  \ee
  where $\oy_{R\eps}$ denotes the solution of the initial value problem
  \bas
	\left\{ \begin{array}{l}
	\oy_{R\eps}'(t)=-c_4 A_{R\eps} \oy_{R\eps}^\frac{p+q}{qa}(t) + c_5(R) A_{R\eps} \eps^\frac{p+q}{a}, 
		\qquad t>0, \\[1mm]
	\oy_{R\eps}(0)=\int_{B_R} u_{0R\eps}^q.
	\end{array} \right.
  \eas
  Now, since our construction of $(u_{0R\eps})_{\eps\in (0,1)}$ guarantees that
  \bas
	A_{R\eps} \to A_R:=\|u_{0R}\|_{L^{q_0}(B_R)}^{-(p+q)(1-a)/(aq_0)}
	\quad \mbox{and} \quad
	\int_{B_R} u_{0R\eps}^q \to \int_{B_R} u_{0R}^q
	\qquad \mbox{as } \eps\searrow 0,
  \eas
  it follows from standard results on continuous dependence that $\oy_{R\eps} \to \oy_R$ in $C^0_{loc}([0,\infty))$ as
  $\eps\searrow 0$, where
  \bas
	\left\{ \begin{array}{l}
	\oy_R'(t)=-c_4 A_R \oy_R^\frac{p+q}{qa}(t), \qquad t>0, \\[1mm]
	\oy_R(0)=\int_{B_R} u_{0R}^q.
	\end{array} \right.
  \eas
  Here an explicit integration shows that if we set $r:=aq/(p+q-aq)$ then
  \bea{1.8}
	\oy_R(t)
	= \bigg( \oy_R^{-\frac{1}{r}}(0) + \frac{c_4}{r} A_R t 
		\bigg)^{-r}
	\le c_6 (A_Rt)^{-r}
	\qquad \mbox{for all } t>0,
  \eea
  where $c_6:=(c_4/r)^{-r}$, and where thanks to (\ref{1.33}),
  \bas
	\frac{1}{r}=\frac{p+q}{qa}-1
	= \frac{2q_0+np}{n(q-q_0)} > 0
  \eas
  according to the fact that $q>q_0$. Therefore, using that $u_{0R} \le u_0$ in estimating
  \bas
	A_R \ge A := \|u_0\|_{L^{q_0}(\R^n)}^{-(p+q)(1-a)/(aq_0)}
	\qquad \mbox{for all } R>0,
  \eas
  letting $\eps\searrow 0$ in (\ref{1.7}) we infer from (\ref{1.8}) that for every $R>0$ we have
  \bas
	\int_{B_R} u_R^r(\cdot,t) \le c_6  (At)^{-(1-\frac{q_0}{q})/(p+\frac{2q_0}{n})}
	\qquad \mbox{for all } t>0,
  \eas
  which on taking the $q$-th root on both sides proves (\ref{1.1}).
\qed
Taking $R\nearrow \infty$ and recalling Lemma~\ref{lem202}, we can thereby already verify the first estimate
in Theorem~\ref{theo200}.
\begin{cor}\label{cor3}
  Let $p\ge 1$ and $u_0\in C^0(\R^n)$ be positive and such that
  \bas
	u_0 \in L^{q_0}(\R^n)
	\qquad \mbox{for some } q_0>0.
  \eas
  Then for any $q>q_0$, the solution $u$ of (\ref{0}) obtained in Proposition~\ref{prop201}
  satisfies $u(\cdot,t)\in L^q(\R^n)$ for all $t>0$, moreover,
  there exists $C(q)>0$ such that
  \be{3.1}
	\|u(\cdot,t)\|_{L^q(\R^n)} \le C(q) t^{-(1-\frac{q_0}{q})/(p+\frac{2q_0}{n})}
	\qquad \mbox{for all } t>0.
  \ee
%
\end{cor}
\proof
  This is an evident consequence of Lemma~\ref{lem1} combined with (\ref{conv_R}) and the monotone convergence theorem.
\qed
\subsection{Upper bounds for $u$ in $L^\infty(\R^n)$. Proof of Theorem~\ref{theo200}}
The constant $C(q)$ inequality gained in Corollary~\ref{cor3} may depend on $q$;
as Proposition~\ref{prop103} will reveal, in general we actually must have $C(q)\to\infty$ as $q\to\infty$.
Accordingly, establishing the respective decay estimates from Theorem~\ref{theo200} which involve the 
norm in $L^\infty(\R^n)$ will, besides Lemma~\ref{lem1}, require an additional ingredient.
This role will be played by the following semi-convexity estimate which is a well-known feature of 
nonlinear diffusion equations of type (\ref{0}). For its derivation, we may thus refer to the literature
(see \cite{aronson_lect} or \cite{win_critexp}, for example).
\begin{lem}\label{lem7}
  Let $R>0$. Then
  \be{7.1}
	\frac{u_{Rt}(x,t)}{u_R(x,t)} \ge - \frac{1}{pt}
	\qquad \mbox{for all $x\in B_R$ and } t>0.
  \ee
\end{lem}
In exploiting this, we shall need the following variant of the Gagliardo-Nirenberg inequality 	
which involves small integrability powers, and in which the dependence of the constant on these powers
is stressed.
\begin{lem}\label{lem6}
  Let $s>2$ be such that $(n-2)s\le 2n$. Then there exists $C(s)>0$ such that for any $\sigma\in (0,2)$ we have
  \bas
	\|\varphi\|_{L^s(\R^n)} \le C(s)\|\nabla \varphi\|_{L^2(\R^n)}^a \|\varphi\|_{L^\sigma(\R^n)}^{1-a}
	\qquad \mbox{for all } \varphi \in W^{1,2}(\R^n) \cap L^\sigma(\R^n),
  \eas
  where
  \bas
	a= a(s,\sigma) :=\left(\frac{n}{\sigma}-\frac{n}{s}\right)\left(1-\frac{n}{2}+\frac{n}{\sigma}\right)^{-1}
	\, \in (0,1].
  \eas
\end{lem}
\proof
  According to the standard Gagliardo-Nirenberg inequality (\cite{friedman}), there exists $c_1(s)>0$ such that
  \be{3.3}
	\|\varphi\|_{L^s(\R^n)} \le c_1(s)\|\nabla \varphi\|_{L^2(\R^n)}^b \|\varphi\|_{L^2(\R^n)}^{1-b}
	\qquad \mbox{for all } \varphi\in C_0^\infty(\R^n),
  \ee
  where 
  \bas
	b=\frac{n}{2}-\frac{n}{s} \in (0,1]
  \eas
  due to our assumptions $s>2$ and $(n-2)s\le 2n$.
  Here we can apply the H\"older inequality to further interpolate
  \be{3.4}
	\|\varphi\|_{L^2(\R^n)}^{1-b} 
	\le \|\varphi\|_{L^s(\R^n)}^{(1-b)d} \|\varphi\|_{L^\sigma(\R^n)}^{(1-b)(1-d)}
	\qquad \mbox{for all } \varphi\in C_0^\infty(\R^n)
  \ee
  with
  \bas
	d= d(s,\sigma):=\left(\frac{1}{\sigma}-\frac{1}{2}\right)\left(\frac{1}{\sigma}-\frac{1}{s}\right)^{-1}.
  \eas
  Combining (\ref{3.3}) with (\ref{3.4}) shows that
  \bas
	\|\varphi\|_{L^s(\R^n)}^{1-(1-b)d} 
	\le c_1(s) \|\nabla \varphi\|_{L^2(\R^n)}^b \|\varphi\|_{L^\sigma(\R^n)}^{(1-b)(1-d)}
	\qquad \mbox{for all } \varphi\in C_0^\infty(\R^n),
  \eas
  so that since
  \bas
	\frac{b}{1-(1-b)d}=a
	\qquad \mbox{and} \qquad
	\frac{(1-b)(1-d)}{1-(1-b)d} = 1-\frac{b}{1-(1-b)d}=1-a,
  \eas
  we obtain
  \bas
	\|\varphi\|_{L^s(\R^n)} \le c_1^\frac{a(s,\sigma)}{b(s)}(s)  \|\nabla \varphi\|_{L^2(\R^n)}^a
	\|\varphi\|_{L^\sigma(\R^n)}^{1-a}
	\qquad \mbox{for all } \varphi\in C_0^\infty(\R^n).
  \eas
\qed
With this tool at hand, we can derive the following consequence of Lemma~\ref{lem7}
which will form the essential step in an iteration procedure to be performed in
Lemma~\ref{lem8}.
\begin{lem}\label{lem9}
  Let $s>2$ be such that $(n-2)s<2n$. Then there exists $C>0$ such that for any $r\ge 1$, the solution $u_R$ of
  (\ref{0R}) satisfies
  \be{9.1}
	\|u_R(\cdot,t)\|_{L^{(p+r)s/2}(B_R)} \le \Big(\frac{Cr}{t} \Big)^\frac{a}{p+r}  
	\|u_R(\cdot,t)\|_{L^r(B_R)}^{1-\frac{pa}{p+r}}
	\qquad \mbox{for all } t>1,
  \ee
  where 
  \be{9.2}
	a= a(r):=\left(\frac{n(p+r)}{2r}-\frac{n}{s}\right)\left(1-\frac{n}{2}+\frac{n(p+r)}{2r}\right)^{-1} \, \in (0,1).
  \ee
\end{lem}
\proof
  We fix a sequence $(\chi_k)_{k\in\N} \subset C^\infty(\R)$ of cut-off functions satisfying $0\le \chi_k \le 1$ and
  $\chi_k' \ge 0$ on $\R$ as well as $\chi_k \equiv 0$ in $(-\infty,1/k)$ and $\chi_k \equiv 1$ in 
  $(2/k,\infty)$.
  Then since $u_R$ is continuous with $u_R|_{\partial B_R}=0$, testing the inequality (\ref{7.1}) by
  $\chi_k(u) u^r$ yields
  \bas
	(p+r-1) \int_{B_R} \chi_k(u_R) u_R^{p+r-2} |\nabla u_R|^2
	+ \int_{B_R} \chi_k'(u_R) u_R^{p+r-1} |\nabla u_R|^2
	\le \frac{1}{pt} \int_{B_R} \chi_k(u_R) u_R^r
	\quad \mbox{for all } t>0.
  \eas
  Since $\chi_k' \ge 0$, $\chi_k \le 1$ and $p+r-1>0$, this implies that
  \bea{9.3}
	\int_{B_R} \chi_k(u_R) u_R^{p+r-2} |\nabla u_R|^2
	&\le& \frac{1}{p(p+r-1)t} \int_{B_R} \chi_k(u_R) u_R^r \nn\\
	&\le& \frac{1}{p(p+r-1)t} \int_{B_R} u_R^r
	\qquad \mbox{for all } t>0.
  \eea
  Thus, if we introduce
  \be{9.33}
	P_k(\xi):=\int_0^\xi \sqrt{\chi_k(\eta)}  \eta^\frac{p+r-2}{2} d\eta, \qquad \xi\ge 0,
  \ee
  and extend $P_k(u_R(\cdot,t))$ by zero so as to become a function $P_k(u_R(\cdot,t)) \in C_0^\infty(\R^n)$, then
  (\ref{9.3}) says that
  \be{9.4}
	\int_{\R^n} |\nabla P_k(u_R)|^2 
	\le \frac{1}{p(p+r-1)t} \int_{B_R} u_R^r
	\qquad \mbox{for all } t>0.
  \ee
  Here we apply the Gagliardo-Nirenberg inequality from Lemma~\ref{lem6} to obtain $c_1>0$ such that for any choice of
  $r>0$, inter alia ensuring that
  \bas
	0<\frac{2r}{p+r}<2,
  \eas
  we have
  \be{9.5}
	\|P_k(u_R)\|_{L^s(\R^n)} \le c_1 \|\nabla P_k(u_R)\|_{L^2(\R^n)}^a\|P_k(u_R)\|_{L^\frac{2r}{p+r}(\R^n)}^{1-a}
	\qquad \mbox{for all } t>0
  \ee
  with $a\in (0,1)$ determined in (\ref{9.2}).
  Again since $0\le \chi_k \le 1$, from (\ref{9.33}) we see that 
  $P_k(\xi) \le \frac{2}{p+r} \xi^\frac{p+r}{2}$ for all $\xi\ge 0$, so that
  \bas
	\|P_k(u_R)\|_{L^\frac{2r}{p+r}(\R^n)}^{1-a}
	\le \bigg\{ \frac{2}{p+r} \Big\| u_R^\frac{p+r}{2} \Big\|_{L^\frac{2r}{p+r}(B_R)} \bigg\}^{1-a} 
	= \Big(\frac{2}{p+r}\Big)^{1-a}  \|u_R\|_{L^r(B_R)}^\frac{(p+r)(1-a)}{2}
  \eas
  for all $t>0$.
  Accordingly, combining (\ref{9.5}) with (\ref{9.4}) shows that
  \bea{9.6}
	\|P_k(u_R)\|_{L^s(\R^n)}
	&\le& \Big(\frac{2}{p+r}\Big)^{1-a}
	 \|\nabla P_K(u_R)\|_{L^2(\R^n)}^a \|u_R\|_{L^r(B_R)}^\frac{(p+r)(1-a)}{2} \nn\\
	&\le& \Big(\frac{2}{p+r}\Big)^{1-a}  \Big(\frac{1}{p(p+r-1)t}\Big)^\frac{a}{2} 
	\bigg(\int_{B_R} u_R^r \bigg)^\frac{a}{2}
	\|u_R\|_{L^r(B_R)}^\frac{(p+r)(1-a)}{2} \nn\\
	&=& \Big(\frac{2}{p+r}\Big)^{1-a}  \Big(\frac{1}{p(p+r-1)t}\Big)^\frac{a}{2} 
	\|u_R\|_{L^r(B_R)}^\frac{p+r-pa}{2}
	\qquad \mbox{for all } t>0.
  \eea
  Now since (\ref{9.33}) along with our choice of $\chi_k$ ensures that for all $\xi<0$ we have
  \bas
	P_k(\xi) \to \frac{2}{p+r} \xi^\frac{p+r}{2}
	\qquad \mbox{as } k\to\infty,
  \eas
  we may apply Fatou's lemma to conclude from (\ref{9.6}) that
  \bas
	\frac{2}{p+r} \|u_R(\cdot,t)\|_{L^\frac{(p+r)s}{2}(B_R)}^\frac{p+r}{2}
	&=& \frac{2}{p+r} \Big\|u_R^\frac{p+r}{2}(\cdot,t)\Big\|_{L^s(B_R)}
	\le \liminf_{k\to\infty} \|P_k(u_R(\cdot,t))\|_{L^s(\R^n)} \\[1mm]
	&\le& \Big(\frac{2}{p+r}\Big)^{1-a}  \Big(\frac{1}{p(p+r-1)t}\Big)^\frac{a}{2}
	 \|u_R(\cdot,t)\|_{L^r(B_R)}^\frac{p+r-pa}{2}
  \eas
  and hence
  \be{9.7}
	\|u_R(\cdot,t)\|_{L^\frac{(p+r)s}{2}(B_R)}
	\le \Big(\frac{p+r}{2}\Big)^\frac{2a}{p+r}  \Big(\frac{1}{p(p+r-1)t}\Big)^\frac{a}{p+r}
	 \|u_R(\cdot,t)\|_{L^r(B_R)}^\frac{p+r-pa}{p+r}
  \ee
  for all $t>0$. Here, we can clearly find $c_2, c_3>1$ such that
  \bas
	\frac{p+r}{2} \le c_2 r
	\quad \mbox{and} \quad
	\frac{1}{p(p+r-1)t} \le \frac{c_3}{rt}
	\qquad \mbox{for all } r\ge 1,
  \eas
  so that (\ref{9.1}) follows from (\ref{9.7}) if we let $C:=c_2^2 c_3$.
\qed
Now an estimate for $u_R$ in $L^\infty(B_R)$ in the flavour of Theorem~\ref{theo200} can be achieved by means of
a Moser-type iteration, making essential use of the dependence of the right-hand side of (\ref{9.1}) on
both $r$ and $t$.
\begin{lem}\label{lem8}
  Let $p\ge 1$ and $u_0\in C^0(\R^n)$ be positive and such that $u_0\in L^{q_0}(\R^n)$ for some $q_0>0$.
 Set $\nu:=\frac{n}{np+2q_0}$.
 Then for any $\delta>0$ there exists $C(\delta)>0$ such that for any $R>0$ we have
  \be{8.1}
	\|u_R(\cdot,t)\|_{L^\infty(B_R)} \le C(\delta)  t^{-\nu+\delta}
	\qquad \mbox{for all } t>0.
  \ee
\end{lem}
\proof
  Let us fix $c_1>0$ such that
  \be{8.2}
	\ln \xi \ge -c_1  (1-\xi)
	\qquad \mbox{for all } \xi\in (1/2,1)
  \ee
  and pick $s>2$ sufficiently close to $2$ such that
$	(n-2)s<2n$.
 Further, given $\delta>0$ satisfying $\delta<\nu$ we can find $\delta' \in (0,\delta)$ such that
 $
	-\nu+\delta
	> - \nu  (1-\delta')
$
  and then take $q>q_0$ large such that still
  \be{8.4}
	-\frac{q-q_0}{q}\nu  (1-\delta')
	< - \nu + \delta.
  \ee
  Next we choose a number $r_0\ge 1$ large enough satisfying
  \be{8.44}
	r_0 \ge q
  \ee
  and
  \be{8.5}
	\frac{p}{p+r_0} < \frac{1}{2}
  \ee
  as well as
  \be{8.6}
	\frac{pc_1}{r_0}  \sum_{i=0}^\infty \Big(\frac{2}{s}\Big)^i
	\le \ln \frac{1}{1-\delta'}
  \ee
  and recursively define
  \bas
	r_{k+1} := \frac{(p+r_k)s}{2} 
  \eas
  for nonnegative integers $k$.
  Then, clearly, $(r_k)_{k\ge 0}$ is increasing with 
  \be{8.7}
	r_k \ge r_0  \Big(\frac{s}{2}\Big)^k
	\qquad \mbox{for all } k\ge 0,
  \ee
  and there exists $c_2>0$ such that
  \be{8.8}
	r_k \le c_2 s^k
	\qquad \mbox{for all } k\ge 0.
  \ee
  In particular, Lemma~\ref{lem9} applies to yield a constant $c_3>1$ such that writing
  \bas
	M_k(t):=\|u_R(\cdot,t)\|_{L^{r_k}(B_R)}, \qquad t>0,
  \eas
  we have
  \be{8.9}
	M_{k+1}(t) \le \Big( \frac{c_3 r_k}{t} \Big)^\frac{a_k}{p+r_k}  M_k^{\theta_k}(t)
	\qquad \mbox{for all $t>0$ and } k\ge 0,
  \ee
  where 
  \be{8.10}
	\theta_k:=1-\frac{pa_k}{p+r_k} \equiv \frac{(1-a_k)p+r_k}{p+r_k}
  \ee
  with
  \be{8.11}
	a_k:=\left(\frac{n(p+r_k)}{2r_k} - \frac{n}{s}\right)\left(1-\frac{n}{2}+\frac{n(p+r_k)}{2r_k}\right)^{-1} \, \in (0,1)
  \ee
  for $k\ge 0$. 
  Now by a straightforward induction, (\ref{8.9}) implies that,
for all $t>0$ and $k\ge 0$, we have
  \be{8.12}
	M_{k+1}(t)	
	\le \Bigg\{ \Big(\frac{c_3}{t}\Big)^{\sum_{j=0}^k \frac{a_j}{p+r_j}  \prod_{i=j+1}^k \theta_i} \Bigg\} 
	\Bigg\{ \prod_{j=0}^k r_j^{\frac{a_j}{p+r_j}  \prod_{i=j+1}^k \theta_i} \Bigg\} 
	M_0^{\prod_{i=0}^k \theta_i}(t),
  \ee
 where we may use that (\ref{8.10}) yields that $\theta_i\in (0,1)$ for all $i\ge 0$ in estimating
  \be{8.122}
	\prod_{i=j+1}^k \theta_i \le 1
	\qquad \mbox{for all $k\ge 0$ and } j\in \{0,\dots,k\}.
  \ee
  Along with (\ref{8.8}), (\ref{8.7}) and the fact that $a_j \le 1$ and $r_j\ge r_0 \ge 1$ for all $j\ge 0$, this yields that
  \bas
	\ln \, \Bigg\{ \prod_{j=0}^k r_j^{\frac{a_j}{p+r_j}  \prod_{i=j+1}^k \theta_i} \Bigg\}
	&\le& \ln \Bigg\{ \prod_{j=0}^k r_j^\frac{a_j}{p+r_j} \Bigg\} 
	= \sum_{j=0}^k \frac{a_j}{p+r_j}  \ln r_j 
	\le \sum_{j=0}^k \frac{1}{r_j}  \ln r_j \\
	&\le& \sum_{j=0}^k \frac{1}{r_0 (\frac{s}{2})^j}  \ln (c_2 s^j) 
	= \frac{\ln c_2}{r_0}  \sum_{j=0}^k \Big(\frac{2}{s}\Big)^j
	+ \frac{\ln s}{r_0}  \sum_{j=0}^k j \Big(\frac{2}{s}\Big)^j
	\qquad \mbox{for all } k\ge 0.
  \eas
  Since both $\sum_{j=0}^\infty (\frac{2}{s})^j$ and $\sum_{j=0}^\infty j (\frac{2}{s})^j$ converge thanks to the 
  fact that $s>2$, from this we infer the existence of $c_4>0$ such that
  \be{8.13}
	\prod_{j=0}^k r_j^{\frac{a_j}{p+r_j}  \prod_{i=j+1}^k \theta_i} \le c_4
	\qquad \mbox{for all } k\ge 0.
  \ee
 Next, (\ref{8.122}) entails that since $c_3>1$ and $\theta_i>0$ for all $i\ge 0$, for every $t\ge 1$ we can estimate
  \bas
	\Big(\frac{c_3}{t}\Big)^{\sum_{j=0}^k \frac{a_j}{p+r_j}  \prod_{i=j+1}^k \theta_i}
	\le c_3^{\sum_{j=0}^k \frac{a_j}{p+r_j}  \prod_{i=j+1}^k \theta_i} 
	\le c_3^{\sum_{j=0}^k \frac{a_j}{p+r_j}}
	\qquad \mbox{for all } k\ge 0,
  \eas
  whence again by using (\ref{8.7}) and that $a_j \le 1$ for all $j\ge 0$ we see that
  \be{8.14}
	\Big(\frac{c_3}{t}\Big)^{\sum_{j=0}^k \frac{a_j}{p+r_j}  \prod_{i=j+1}^k \theta_i}
	\le c_3^{\frac{1}{r_0}  \sum_{j=0}^k (\frac{2}{s})^j}
	\le c_5:=c_3^{\frac{1}{r_0} \sum_{j=0}^\infty (\frac{2}{s})^j}
	\qquad \mbox{for all $t\ge 1$ and } k\ge 0.
  \ee
  Finally, to control the rightmost factor in (\ref{8.12}) we first apply
  Lemma~\ref{lem1} to find $c_6>1$ such that
  \bas
	M_0(t) \le c_6 t^{-\frac{r_0-q_0}{r_0}\nu}
	\qquad \mbox{for all } t>0,
  \eas
  which according to (\ref{8.44}) implies that
  \be{8.144}
	M_0(t) \le c_6 t^{-\frac{q-q_0}{q}\nu}
	\qquad \mbox{for all } t\ge 1.
  \ee
  Now whenever $t\ge 1$ is such that the right-hand side herein is larger than $1$, that is, when
  \bas
	1 \le t < t_0:=c_6^\frac{q}{\nu(q-q_0)},
  \eas
  from (\ref{8.12}), (\ref{8.13}) and (\ref{8.14}) we trivially infer that
  \be{8.15}
	M_{k+1}(t) \le c_4 c_5 c_6 \le c_4 c_5 c_6 t_0^{\nu-\delta}
	 t^{-\nu + \delta},
  \ee
  again because $\theta_i\in (0,1)$ for all $i\ge 0$. 
  If, conversely, $t\ge t_0$ then we estimate
  \be{8.16}
	M_0^{\prod_{i=0}^k \theta_i} (t) \le \left( c_6 t^{-\frac{q-q_0}{q}\nu} 
		\right)^{\prod_{i=0}^k \theta_i},
  \ee
  where, since $\theta_i\ge 1-\frac{p}{p+r_i}$ for all $i\ge 0$, we have
  \bas
	\theta_i \ge 1-\frac{p}{p+r_0}>\frac{1}{2}
	\qquad \mbox{for all } i\ge 0
  \eas
  according to (\ref{8.11}), (\ref{8.7}) and (\ref{8.5}), using (\ref{8.2}) and once more (\ref{8.7}) we know that
  \bas
	\ln \Bigg\{ \prod_{i=0}^k \theta_i \Bigg\}
	&=& \sum_{i=0}^k \ln \theta_i 
	\ge -c_1  \sum_{i=0}^k (1-\theta_i) 
	\ge - c_1  \sum_{i=0}^k \frac{p}{p+r_i} 
	\ge -pc_1  \sum_{i=0}^k \frac{1}{r_i} \\
	&\ge& - \frac{pc_1}{r_0}  \sum_{i=0}^k \Big(\frac{2}{s}\Big)^i 
	\ge -\ln \frac{1}{1-\delta'}
	\qquad \mbox{for all } k\ge 0
  \eas
  as a consequence of (\ref{8.6}). Therefore,
  \bas
	\prod_{i=0}^k \theta_i \ge 1-\delta'
	\qquad \mbox{for all } k\ge 0,
  \eas
  so that (\ref{8.16}) combined with (\ref{8.4}) and the fact that $t_0>1$ guarantees that
  \bas
	M_0^{\prod_{i=0}^k \theta_i} (t)
	\le \left( c_6 t^{-\frac{q-q_0}{q}\nu} \right)^{1-\delta'} 
	\le c_6^{1-\delta'}  t^{-\nu +\delta}
	\qquad \mbox{for all } t\ge t_0.
  \eas
  In conjunction with (\ref{8.15}) this establishes (\ref{8.1}) upon an evident choice of $C(\delta)$.
\qed
In the limit $R\nearrow \infty$, this proves the estimate (\ref{200.2}) in
Theorem~\ref{theo200}.
\begin{cor}\label{cor10}
  Assume that $p\ge 1$, and that $u_0\in C^0(\R^n)$ is positive fulfilling $u_0\in L^{q_0}(\R^n)$ for some $q_0>0$.
  Then for any $\delta>0$ one can find $C(\delta)>0$ such that the solution $u$ of (\ref{0}) 
  obtained in Proposition~\ref{prop201} satisfies
  \bas
	\|u(\cdot,t)\|_{L^\infty(\R^n)} \le C(\delta)  t^{-\nu+\delta}
	\qquad \mbox{for all } t \ge 1.
  \eas
\end{cor}
\proof
  The claim immediately results from Lemma~\ref{lem8} upon taking $R\nearrow \infty$.
\qed
On furthermore taking $q_0\searrow 0$ herein, we can also deduce the estimate (\ref{200.3}) on fast decay in
$L^\infty(\R^n)$ for rapidly decreasing initial data.
\begin{cor}\label{cor11}
  Let $p\ge 1$ and $u_0\in C^0(\R^n)$ be positive and such that
$u_0\in L^q(\R^n)$ for all $q>0$.
  Then for any $\delta>0$ there exists $C(\delta)>0$ such that for the solution $u$ of (\ref{0}) from 
Proposition~\ref{prop201} we have
  \be{11.1}
	\|u(\cdot,t)\|_{L^\infty(\R^n)}
	\le C(\delta)  t^{-\frac{1}{p}+\delta}
	\qquad \mbox{for all } t\ge 1.
  \ee
\end{cor}
\proof
  Given $\delta>0$, we fix $q_0>0$ small enough such that
  $-\nu < -\frac{1}{p}+\delta$,
  and then choose $\delta'>0$ fulfilling
$-\nu + \delta' \le -\frac{1}{p}+\delta$.
  An application of Corollary~\ref{cor10} 		
  then yields (\ref{11.1}).
\qed
Summarizing, we can thereby complete the proof of Theorem~\ref{theo200}.\abs
\proofc of Theorem~\ref{theo200}.\quad
  We only need to collect the results provided by Corollary~\ref{cor3},
Corollary~\ref{cor10} and Corollary~\ref{cor11}.
\qed
\mysection{Estimates from below. Optimality}
In order to derive lower estimates for solutions of (\ref{0}), given any such solution $u$,
for $x\in\R^n$ and $t\ge 0$ we introduce
\be{transf}
	v(x,\tau):=(t+1)^\frac{1}{p}  u(x,t),
	\qquad \tau=\ln (t+1),
\ee
whence, as can easily be verified, $v$ is a positive classical solution of
\be{0v}
	\left\{ \begin{array}{ll}
	v_\tau=v^p \Delta v + \frac{1}{p} v, \qquad & x\in \R^n, \ t>0, \\[1mm]
	v(x,0)=u_0(x), \qquad & x\in\R^n.
	\end{array} \right.
\ee
The following result can be obtained by a comparison argument from below with functions of the form 
$\uv(x,\tau):=y_R(\tau)  (\Theta_R(x)+\delta_R)$, where $\Theta_R$ denotes the principal eigenfunction of the 
Dirichlet Laplacian in $B_R$, and where the number $\delta_R>0$ and the positive function $y_R$ are chosen appropriately
(cf.~\cite[Lemma~2.1]{win_critexp}).
\begin{lem}\label{lem102}
  Let $p\ge 1$, and suppose that $u_0\in C^0(\R^n)$ is positive and that $u$ is a globally defined positive
  classical solution of (\ref{0}). Then for every $R>0$,
  the function $v$ defined by (\ref{transf}) satisfies
  \bas
	\inf_{x\in B_R} v(x,t) \to + \infty
	\qquad \mbox{as } t\to\infty.
  \eas
\end{lem}
\proofc of Proposition~\ref{prop103}.\quad
  The claimed divergence property is an immediate consequence of Lemma~\ref{lem102} and (\ref{transf}).
\qed
\subsection{Initial data with algebraic decay}
Let us next focus on the particular case when the initial data essentially decay algebraically in space.
The key step in our analysis of the corresponding solutions of (\ref{0}) will consist in a comparison
argument involving separated solutions of the PDE in (\ref{0}) in suitably chosen balls. 
The spatial profiles arising therein are addressed in the following lemma asserting a favorable
scaling property.
\begin{lem}\label{lem12}
  Let $p\ge 1$. For $R>0$, let $w_R \in C^0(\bar B_R) \cap C^2(B_R)$ denote the positive solution of
  \be{wR}
	\left\{ \begin{array}{ll}
	-\Delta w_R = \frac{1}{p} w_R^{1-p}, \qquad &x\in B_R, \\[1mm]
	w_R=0, \qquad &x\in \partial B_R.
	\end{array} \right.
  \ee
  Then for every $R>0$ we have
  \be{12.1}
	w_R(x)=R^\frac{2}{p}  w_1\Big(\frac{x}{R}\Big)
	\qquad \mbox{for all } x\in B_R.
  \ee
\end{lem}
\proof 
  Since 
  \bas
	\Delta \Big\{ R^\frac{2}{p}  w_1 \Big(\frac{x}{R}\Big)\Big\}
	+ \frac{1}{p}  \Big\{ R^\frac{2}{p}  w_1\Big(\frac{x}{R}\Big)\Big\}
	&=& R^{\frac{2}{p}-2}  (\Delta w_1)\Big(\frac{x}{R}\Big)
	+ \frac{1}{p} R^\frac{2(1-p)}{p} w_1^{1-p}\Big(\frac{x}{R}\Big) \\
	&=& R^{\frac{2}{p}-2}  \Big\{ \Delta w_1 + \frac{1}{p} w_1^{1-p} \Big\} \Big(\frac{x}{R}\Big)
	\qquad \mbox{for all } x\in B_R,
  \eas
  it results from the properties defining $w_1$ that $B_R \ni x \mapsto R^\frac{2}{p} w_1(\frac{x}{R})$ solves (\ref{wR}).
  By uniqueness of positive solutions to (\ref{wR}) (\cite{wiegner}), (\ref{12.1}) thus follows.
\qed
Now by an appropriate comparison argument from below we shall achieve the following quantitative 
implication of an algebraic lower decay estimate on $u_0$ for the size of $v$ in balls with a certain 
time-dependent radius.
\begin{lem}\label{lem13}
  Let $p\ge 1$, and suppose that $u_0\in C^0(\R^n)$ is such that there exist $\gamma>0$ and $c_0>0$ fulfilling
  \be{13.1}
	u_0(x) \ge C_0  (1+|x|)^{-\gamma}
	\qquad \mbox{for all } x\in \R^n.
  \ee
  Then there exists $C>0$ such that for $v$ as given by (\ref{transf}) we have
  \be{13.2}
	v(x,\tau) \ge C  w_{R(\tau)} (x)
	\qquad \mbox{for all $x\in B_{R(\tau)}$ and } \tau>0,
  \ee
  where $w_{R(\tau)}$ denotes the positive solution of (\ref{wR}) corresponding to
  $R(\tau):=e^{\tau/(p\gamma+2)}$.
\end{lem}
\proof
  We fix an arbitrary $\tau_0>0$ and let
  \be{13.4}
	R:=e^\frac{\tau_0}{p\gamma+2}.
  \ee
  Then writing $c_1:=\|w_1\|_{L^\infty(B_1)}$, we set
  \be{13.5}
	\delta:=\frac{C_0}{2^\gamma c_1}  R^{-\gamma-\frac{2}{p}}
  \ee
  and define $y\in C^1([0,\infty))$ to be the solution of
  \be{13.6}
	\left\{ \begin{array}{l}
	y'(\tau)=\frac{1}{p} y(\tau) - \frac{1}{p} y^{p+1}(\tau), \qquad \tau>0, \\[1mm]
	y(0)=\delta,
	\end{array} \right.
  \ee
  that is, we let
  \be{13.7}
	y(\tau):= \Big\{ \delta^{-p}  e^{-\tau} + 1 - e^{-\tau}\Big\}^{-\frac{1}{p}}
	\qquad \mbox{for } \tau\ge 0.
  \ee
  We finally introduce
  \bas
	\uv(x,\tau):=y(\tau)  w_R(x),
	\qquad x\in \bar B_R, \ \tau\ge 0,
  \eas
  and first observe that clearly $\uv = 0 < v$ on $\partial B_R \times (0,\infty)$.  
  Moreover, since from Lemma~\ref{lem12} we know that
  \bas
	w_R(x) \le R^\frac{2}{p} \|w_1\|_{L^\infty(B_1)} = c_1 R^\frac{2}{p}
	\qquad \mbox{for all } x\in B_R,
  \eas
  and since (\ref{13.1}) along with the fact that $R\ge 1$ implies that
  \bas
	u_0(x) \ge C_0  (1+|x|)^{-\gamma} \ge 2^{-\gamma} C_0 R^{-\gamma}
	\qquad \mbox{for all } x\in B_R,
  \eas
  thanks to (\ref{13.5}) we obtain that
  \bas
	\frac{\uv(x,0)}{v(x,0)}
	= \frac{\delta w_R(x)}{u_0(x)} 
	\le \frac {2^\gamma c_1}{C_0} R^{\gamma+\frac{2}{p}} \delta 
	= 1
	\qquad \mbox{for all } x\in B_R.
  \eas
  As (\ref{13.6}) entails that furthermore
  \bas
	\uv_t - \uv^p \Delta \uv - \frac{1}{p} \uv
	= y' w_R - y^{p+1} w_R^p \Delta w_R - \frac{1}{p} y w_R 
	= \Big\{ y' + \frac{1}{p}y^{p+1} - \frac{1}{p} y \Big\}  w_R
	= 0
	\quad \mbox{in } B_R \times (0,\infty),
  \eas
  the comparison principle (cf.~\cite{wiegner} for a version adequate for the present purpose)
  guarantees that
  \be{13.8}
	\uv \le v
	\qquad \mbox{in } B_R\times (0,\infty).
  \ee
  Since from (\ref{13.7}) we see that
  \bas
	y(\tau) \ge \Big\{ \delta^{-p} e^{-\tau} + 1 \Big\}^{-\frac{1}{p}}
	\qquad \mbox{for all } \tau>0,
  \eas
  and that hence, by (\ref{13.4}) and (\ref{13.5}), 
  \bas
	y(\tau_0) &\ge& \Big\{\delta^{-p} e^{-\tau_0}+1\Big\}^{-\frac{1}{p}} 
	= \bigg\{ \Big(\frac{C_0}{2^\gamma c_1}\Big)^{-p} R^{p\gamma+2}  R^{-(p\gamma+2)} + 1 \bigg\}^{-\frac{1}{p}} 
	= \bigg\{ \Big(\frac{2^\gamma c_1}{C_0}\Big)^p + 1 \bigg\}^{-\frac{1}{p}},
  \eas
  evaluating (\ref{13.8}) at $\tau=\tau_0$ readily establishes (\ref{13.2}).
\qed
Properly evaluating this lower estimate we obtain the claimed estimate from below for solutions which initially 
decay no faster than algebraically.\abs
%
%
\proofc of Theorem~\ref{theo2000}~(i).\quad 
  According to Lemma~\ref{lem13} and (\ref{transf}), (\ref{2000.3}) implies that
  there exists  $c_1>0$ such that 
  \be{104.5}
	u(x,t) \ge c_1 (t+1)^{-\frac{1}{p}}  w_{R(\tau)}(x)
	\qquad \mbox{for all $x\in B_{R(\tau)}$ and } t>0,
  \ee
  where $R(\tau)=e^\frac{\tau}{p\gamma+2}$ with $\tau=\ln(t+1)$ and $w_{R(\tau)}$ is taken from (\ref{wR}).
  As for the case $q=\infty$, in light of Lemma~\ref{lem12} this immediately implies that
  \bas
	\|u(\cdot,t)\|_{L^\infty(\R^n)}
	\ge u(0,t) 
	\ge c_1 (t+1)^{-\frac{1}{p}}  R^\frac{2}{p}(\tau) w_1(0) 
	= c_1 w_1(0) (t+1)^{-\frac{\gamma}{p\gamma+2}}\qquad \mbox{for all } t\ge 0
  \eas
and thereby proves (\ref{2000.4}) in this particular situation.\\
  For finite $q$, an integration of (\ref{104.5}) using the substitution $y=\frac{x}{R}$ shows that
  \bas
	\int_{\R^n} u^q(x,t)dx
	&\ge& c_1^q (t+1)^{-\frac{q}{p}} \int_{B_{R(\tau)}} w^q_{R(\tau)} (x) dx \\
	&=& c_1^q (t+1)^{-\frac{q}{p}} 
	\int_{B_{R(\tau)}} \bigg\{ R^\frac{2}{p}(\tau)  w_1 \Big(\frac{x}{R}\Big) \bigg\}^q dx \\
	&=& c_1^q (t+1)^{-\frac{q}{p}}  R^{\frac{2q}{p}+n}(\tau)  \int_{B_1} w_1(y) dy \\
	&=& c_1^q (t+1)^{-\frac{q\gamma-n}{p\gamma+2}} 
	\qquad \mbox{for all } t>0.
  \eas
  Taking the $q$-th root here we readily arrive at (\ref{2000.4}) also in this case.
\qed
%
%
%
%
%
%
%
%
%
%
%
%
%
%
\mysection{Optimality of $L^q$ decay estimates. Proof of Theorem~\ref{theo100}}
Another application of Theorem \ref{theo2000} (i) to suitably chosen $\gamma$ finally shows optimality
of the estimates (\ref{200.1}) and (\ref{200.2}) as formulated in Theorem~\ref{theo100}.\abs
\proofc of Theorem~\ref{theo100}.\quad
  Since
  \bas
	\frac{q\gamma-n}{q(p\gamma+2)}\quad \to \quad
	\frac{n(q-q_0)}{q(np+2q_0)}
	\qquad \mbox{as } \gamma \searrow \frac{n}{q_0},
  \eas
  we see that given $q_0>0$, $q \in (q_0,\infty]$ and $\delta>0$ we can fix
  $\gamma>n/q_0$ such that
  \bas
	\frac{q\gamma-n}{q(p\gamma+2)}<\gamma(q_0,q)+\delta.
  \eas
  We then choose any positive function $u_0\in C^0(\R^n)$ with the property that 
  \be{100.4}
	c_1  (1+|x|)^{-\gamma} \le u_0(x) \le c_2  (1+|x|)^{-\gamma}
	\qquad \mbox{for all } x\in\R^n
  \ee
  with certain positive constants $c_1$ and $c_2$.
  Since $\gamma>\frac{n}{q_0}$, the right inequality herein ensures that $u_0\in L^{q_0}(\R^n)$,
  whereas the left allows for an application of Theorem~\ref{theo2000} (i) which provides $c_3>0$ fulfilling
  \bas
	\|u(\cdot,t)\|_{L^q(\Omega)} \ge c_3 t^{-\frac{\gamma-\frac{n}{q}}{p\gamma+2}}
	\qquad \mbox{for all } t>1.
  \eas
  Thanks to (\ref{100.4}), this proves (\ref{100.1}).
\qed
\mysection{A family of self-similar solutions of (\ref{0}) with algebraic decay}
In the following, given $p\ge 1, \alpha>0, \beta>0$ and $A>0$ we let $f\equiv f_A\equiv \f$ denote the solution
of the initial value problem (\ref{0f})
which we suppose to be defined as a positive classical solution in $[0,\xi_0)$ with maximally chosen $\xi_0\in (0,\infty]$.
\begin{lem}\label{lem51}
  Let $p\ge 1$, $\alpha>0$ and $\beta>0$. Then for any choice of $A>0$, $f_A$ exists on all of $[0,\infty)$ and satisfies
  \be{51.1}
	f_A(\xi)>0
	\quad \mbox{and} \quad
	f_A'(\xi) \le 0
	\qquad \mbox{for all } \xi\ge 0.
  \ee
\end{lem}
\proof
  Since $\alpha>0$, it is clear that $f:=f_A$ cannot attain a positive local minimum on $(0,\infty)$
  and hence satisfies $f>0$ as well as $f'\le 0$ on $[0,\xi_0)$, whence we only need to show that $\xi_0=\infty$.
  Thus assuming for contradiction that $\xi_0$ be finite, we accordingly obtain that $f(\xi)\searrow 0$ as 
  $\xi\nearrow \xi_0$, and for $\xi\in (0,\xi_0)$ we may divide the ODE in (\ref{0f}) by
  $\xi f(\xi)$
  to see that
  \bas
	\beta (\ln f)' = -\frac{\alpha}{\xi} - \xi^{-n} f^{p-1}  (\xi^{n-1} f')'
	\qquad \mbox{for all } \xi\in (0,\xi_0),
  \eas
  and that hence with $\xi_1:=\xi_0/2$ we have
  \be{51.11}
	\beta \ln f(\xi)
	= \beta\ln f(\xi_1)
	- \alpha \ln \frac{\xi}{\xi_1}
	- \int_{\xi_1}^\xi \sigma^{-n} f^{p-1}(\sigma)  (\sigma^{n-1} f'(\sigma))' d\sigma
	\qquad \mbox{for all } \xi\in (\xi_1,\xi_0).
  \ee
  Here an integration by parts shows that
  \bea{51.2}
	- \int_{\xi_1}^\xi \sigma^{-n} f^{p-1}(\sigma)  (\sigma^{n-1} f'(\sigma))' d\sigma
	&=& (p-1) \int_{\xi_1}^\xi \frac{1}{\sigma} f^{p-2}(\sigma) f'^2(\sigma) d\sigma 
	- n\int_{\xi_1}^\xi \frac{1}{\sigma^2} f^{p-1}(\sigma) f'(\sigma) d\sigma \nn\\
	- \frac{1}{\xi} f^{p-1}(\xi) f'(\xi)
	+ \frac{1}{\xi_1} f^{p-1}(\xi_1) f'(\xi_1)
	&\ge&- n\int_{\xi_1}^\xi \frac{1}{\sigma^2} f^{p-1}(\sigma) f'(\sigma) d\sigma 
	-c_1
	\qquad \mbox{for all } \xi\in (\xi_1,\xi_0)
  \eea
  with $c_1:= f^{p-1}(\xi_1) f'(\xi_1)/\xi_1 \ge 0$, because $p\ge 1$ and $f' \le 0$ on $(0,\xi_0)$.
  Once more integrating by parts, since $f(\xi)\to 0$ as $\xi\nearrow\xi_0$ we obtain that
  \bas
	&&n\int_{\xi_1}^\xi \frac{1}{\sigma^2} f^{p-1}(\sigma) f'(\sigma) d\sigma 
	= \frac{n}{p} \int_{\xi_1}^\xi \frac{1}{\sigma^2}  (f^p)'(\sigma) d\sigma 
	= \frac{2n}{p} \int_{\xi_1}^\xi \frac{1}{\sigma^3} f^p(\sigma) d\sigma
	+ \frac{n}{p\xi^2} f^p(\xi)
	- \frac{n}{p\xi_1^2} f^p(\xi_1) \\
	&&\le \frac{2n}{p} \int_{\xi_1}^\xi \frac{1}{\sigma^3} f^p(\sigma) d\sigma
	+ \frac{n}{p\xi^2} f^p(\xi) \quad 
	\to \quad  \frac{2n}{p} \int_{\xi_1}^{\xi_0} \frac{1}{\sigma^3} f^p(\sigma) d\sigma 
	=: c_2
	\qquad \mbox{as } \xi\nearrow\xi_0,
  \eas
  where the fact that $\xi_1\in (0,\xi_0)$ clearly asserts that $c_2$ is finite. 
  From (\ref{51.11}) and (\ref{51.2}) we therefore conclude that
  \bas
	\liminf_{\xi\nearrow\xi_0} \Big\{ \beta \ln f(\xi) \Big\}
	\ge \beta \ln f(\xi_1) - \alpha\ln \frac{\xi_0}{\xi_1} - c_1 - c_2,
  \eas
  which together with our assumption $\beta>0$ contradicts the hypothesis that $f(\xi)\to 0$ as $\xi\nearrow \xi_0$
  and thereby proves that actually $\xi_0=\infty$, as claimed.
\qed
In the proofs of Lemma~\ref{lem53} and Lemma~\ref{lem56} we shall use the
following:
\begin{lem}\label{lem57}
  Let $p>1$, $\alpha>0$, $\beta>0$ and $A>0$. Then the solution $f_A$ of (\ref{0f}) satisfies
  \be{57.1}
	\xi^{n-1} f'(\xi)
	+ \frac{\beta}{p-1} \Big(n+\frac{(p-1)\alpha}{\beta}\Big) \int_0^\xi \sigma^{n-1} f^{1-p}(\sigma)d\sigma 
	- \frac{\beta}{p-1} \xi^n f^{1-p}(\xi)
	=0
	\qquad \mbox{for all } \xi>0.
  \ee
\end{lem}
\proof
  We again let $f:=f_A$ and rewrite (\ref{0f}) in the form
  \bas
	0 = (\xi^{n-1}f')' + \beta \xi^n f^{-p} f' + \alpha \xi^{n-1} f^{1-p} 
	= (\xi^{n-1} f')' - \frac{\beta}{p-1} \xi^{n-1}  
	\bigg\{ \xi (f^{1-p})' - \frac{(p-1)\alpha}{\beta} f^{1-p} \bigg\},
  \eas
  whence computing
  \bas
	\xi^{\frac{(p-1)\alpha}{\beta}+1} \Big( \xi^{-\frac{(p-1)\alpha}{\beta}} f^{1-p} \Big)'
	= \xi (f^{1-p})' - \frac{(p-1)\alpha}{\beta} f^{1-p},
  \eas
  we obtain that
  \bas
	0=(\xi^{n-1} f')'
	- \frac{\beta}{p-1} \xi^{n+\frac{(p-1)\alpha}{\beta}}  \Big( \xi^{-\frac{(p-1)\alpha}{\beta}} f^{1-p} \Big)'
	\qquad \mbox{for all } \xi>0.
  \eas
  Using that $f'(0)=0$, by integration we achieve the identity
  \be{53.3}
	0 = \xi^{n-1} f'(\xi)
	- \frac{\beta}{p-1} \int_0^\xi \sigma^{n+\frac{(p-1)\alpha}{\beta}} 
	 \Big(\sigma^{-\frac{(p-1)\alpha}{\beta}} f^{1-p}\Big)'(\sigma) d\sigma
	\qquad \mbox{for all } \xi>0,
  \ee
  where integrating by parts we find that
  \bas
	-\frac{\beta}{p-1} \int_0^\xi \sigma^{n+\frac{(p-1)\alpha}{\beta}}
	 \Big( \sigma^{-\frac{(p-1)\alpha}{\beta}} f^{1-p}\Big)'(\sigma) 
	&=& \frac{\beta}{p-1} \Big(n+\frac{(p-1)\alpha}{\beta}\Big) \int_0^\xi \sigma^{n-1} f^{1-p}(\sigma)d\sigma \\
	& & - \frac{\beta}{p-1} \xi^n f^{1-p}(\xi)
	\qquad \mbox{for all } \xi>0,
  \eas
  so that (\ref{57.1}) is equivalent to (\ref{53.3}).
\qed
\subsection{A lower bound for $\f$}
A preliminary estimate for $f_A$ from below near $\xi=0$ can be obtained by quite straightforward
manipulations of (\ref{0f}).
\begin{lem}\label{lem52}
  Let $p\ge 1$, $\alpha>0$ and $\beta>0$. 
  Then there exist $\xis>0$ and $K>0$ such that for any $A>0$, the solution $f_A$ of (\ref{0f}) satisfies
  \be{52.11}
	f_A(\xi) \ge \frac{A}{2}
	\qquad \mbox{for all } \xi\in [0,\xis],
  \ee
  and that moreover
  \be{g}
	g_A(\xi):=\int_0^\xi \sigma^{n-1} f_A^{1-p}(\sigma)d\sigma,
	\qquad \xi\ge 0, \ A>0,
  \ee
  has the property that
  \be{52.1}
	\xis^{-n-\frac{(p-1)\alpha}{\beta}} g_A(\xis) \le K A^{1-p}
	\qquad \mbox{for all } A>0.
  \ee
\end{lem}
\proof
  Using that $f:=f_A$ is nonincreasing by Lemma~\ref{lem51}, from (\ref{0f}) we obtain that
  \bas
	f''=-\frac{n-1}{\xi} f' - \beta\xi f' - \alpha f \ge -\alpha f
	\qquad \mbox{for all } \xi>0,
  \eas
  and that hence $\frac{1}{2} (f'^2)' \ge -\frac{\alpha}{2} (f^2)'$ on $(0,\infty)$.
  On integration, this yields
  \bas
	\frac{1}{2} f'^2(\xi) \ge -\frac{\alpha}{2} f^2(\xi) + \frac{\alpha}{2} A^2
	\qquad \mbox{for all } \xi>0
  \eas
  and thus
  \bas
	f'(\xi) \ge -\sqrt{\alpha} \sqrt{A^2-f^2(\xi)}
	\qquad \mbox{for all } \xi>0.
  \eas
  By an ODE comparison, we thereby conclude that
  \bas
	f(\xi) \ge A \, \cos (\sqrt{\alpha}\xi)
	\qquad \mbox{for all } \xi\in \Big(0,\frac{\pi}{2\sqrt{\alpha}}\Big),
  \eas
  which means that if we let $\xis:=\frac{\pi}{3\sqrt{\alpha}}$, then
  \bas
	f(\xi)\ge A\, \cos \frac{\pi}{3} = \frac{A}{2}
	\qquad \mbox{for all }  \xi\in (0,\xis).
  \eas
  This establishes (\ref{52.11}), and since $p\ge 1$ this moreover entails that
  \bas
	\int_0^{\xis} \sigma^{n-1} f^{1-p}(\sigma) d\sigma
	\le \frac{2^{p-1} \xis^n}{n}  A^{1-p},
  \eas
  which readily proves (\ref{52.1}).
\qed
Making use of this information, on the basis of (\ref{57.1}) we can make sure that in the case $p>1$,
a lower bound for the asymptotic behaviour of $f_A$ is determined by the decay of positive solutions to
the first-order equation $\beta \xi f' + \alpha f=0$:
\begin{lem}\label{lem53}
  Let $p>1$, and suppose that $\alpha>0$ and $\beta>0$.
  Then there exists $L>0$ such that for every $A>0$, the solution $f_A$ of (\ref{0f}) satisfies
  \be{53.1}
	f_A(\xi) \ge LA (1+\xi)^{-\frac{\alpha}{\beta}}
	\qquad \mbox{for all } \xi\ge 0.
  \ee
\end{lem}
\proof
  Reformulating (\ref{57.1}) in terms of $g(\xi):=\int_0^\xi \sigma^{n-1} f^{1-p}(\sigma)d\sigma$, $\xi\ge 0$, 
  since $f'\le 0$ on $[0,\infty)$ we obtain that
  \be{53.5}
	\Big(n+\frac{(p-1)\alpha}{\beta}\Big) g(\xi) - \xi^n f^{1-p}(\xi) \ge 0
	\qquad \mbox{for all } \xi>0
  \ee
  and that hence
  \bas
	\xi g'(\xi) \le \Big(n+\frac{(p-1)\alpha}{\beta}\Big) g(\xi)
	\qquad \mbox{for all } \xi>0.
  \eas
  By e.g.~a comparison argument, from this we infer that with $\xis>0$ and $K>0$ as provided by Lemma~\ref{lem52} we have
  \bas
	g(\xi) \le g(\xis)  \Big(\frac{\xi}{\xis}\Big)^{n+\frac{(p-1)\alpha}{\beta}} 
	\le K A^{1-p}  \xi^{n+\frac{(p-1)\alpha}{\beta}}
	\qquad \mbox{for all } \xi>\xis.
  \eas
  Going back to (\ref{53.5}), we see that this yields information on $f$ itself by implying the inequality
  \bas
	\xi^n f^{1-p}(\xi)
	\le \Big(n+\frac{(p-1)\alpha}{\beta}\Big) g(\xi) 
	\le \Big(n+\frac{(p-1)\alpha}{\beta}\Big) K A^{1-p}
	\xi^{n+(p-1)\alpha/\beta}
	\qquad \mbox{for all } \xi>\xis,
  \eas
  which is equivalent to
  \bas
	f(\xi) \ge \Big(n+\frac{(p-1)\alpha}{\beta}\Big)^{-\frac{1}{p-1}} K^{-\frac{1}{p-1}} 
	\xi^{-\frac{\alpha}{\beta}}
	\qquad \mbox{for all } \xi>\xis.
  \eas
  This shows that if $L>0$ is sufficiently small then the estimate in (\ref{53.1}) holds for all $\xi>\xis$,
  whereas (\ref{52.11}) guarantees that on diminishing $L$ if necessary we can achieve the desired lower bound also for
  $\xi\in [0,\xis]$.
\qed
\subsection{An upper bound for $\f$}
Again by means of (\ref{57.1}), we can complete our description of the asymptotic decay of $f_A$ for 
$p>1$ as follows.
\begin{lem}\label{lem56}
  Let $p>1, \alpha>0, \beta>0$ and $A>0$.
  Then there exists $C>0$ with the property that the solution $f_A$ of (\ref{0f}) satisfies
  \be{56.1}
	f_A(\xi) \le C  (1+\xi)^{-\frac{\alpha}{\beta}}
	\qquad \mbox{for all } \xi\ge 0.
  \ee
\end{lem}
\proof
  Again abbreviating $f:=f_A$ and $g(\xi):=\int_0^\xi \sigma^{n-1} f^{1-p}(\sigma)d\sigma$ for $\xi\ge 0$,
  from (\ref{57.1}) we obtain
  \bas
	c_1 \xi^{n-1} f'(\xi)
	+ \Big(n+\frac{(p-1)\alpha}{\beta}\Big) g(\xi) - \xi g'(\xi)=0
	\qquad \mbox{for all } \xi>0
  \eas
  with $c_1:=\frac{p-1}{\beta}$, that is,
  \be{56.2}
	g'(\xi)=\Big(n+\frac{(p-1)\alpha}{\beta}\Big)  \frac{g(\xi)}{\xi} + c_1 \xi^{n-2} f'(\xi)
	\qquad \mbox{for all } \xi>0.
  \ee
  In order to prepare an integration thereof, let us first make sure that if we define
  \be{56.3}
	\xi_1:=\sqrt{2nc_1 A^p},
  \ee
  then
  \be{56.4}
	c_1 \xi_1^{-2} f(\xi_1) \le \frac{1}{2} \xi_1^{-n} g(\xi_1).
  \ee
  Indeed, this follows from the observation that using $f(\sigma)\le A$ for all $\sigma\ge 0$ by
  Lemma~\ref{lem51}, 
  we can estimate
  \bas
	2c_1 \xi_1^{n-2} \frac{f(\xi_1)}{g(\xi_1)} 
	= 2c_1 \xi_1^{n-2}  \frac{f(\xi_1)}{\int_0^{\xi_1} \sigma^{n-1} f^{1-p}(\sigma)d\sigma} 
	\le 2c_1 \xi_1^{n-2}  \frac{A}{A^{1-p} \int_0^{\xi_1} \sigma^{n-1}d\sigma}
	= \frac{2nc_1 A^p}{\xi_1^2} 
	= 1
  \eas
  thanks to (\ref{56.3}).\\
  We now invoke the variation-of-constants formula associated with the initial value problem for (\ref{56.2}) starting
  from $\xi=\xi_1$ to see that
  \bea{56.5}
	g(\xi)
	&=& g(\xi_1)  e^{\left(n+\frac{(p-1)\alpha}{\beta}\right) \int_{\xi_1}^\xi \frac{d\sigma}{\sigma}}
	+ c_1 \int_{\xi_1}^\xi e^{\left(n+\frac{(p-1)\alpha}{\beta}\right) \int_\sigma^\xi \frac{d\tau}{\tau}}
	 \sigma^{n-2} f'(\sigma) d\sigma \nn\\
	&=& g(\xi_1)  \Big(\frac{\xi}{\xi_1} \Big)^{n+\frac{(p-1)\alpha}{\beta}}
	+ c_1 \xi^{n+\frac{(p-1)\alpha}{\beta}} \int_{\xi_1}^\xi \sigma^{-\frac{(p-1)\alpha}{\beta}-2} f'(\sigma) d\sigma
	\qquad \mbox{for all } \xi>\xi_1,
  \eea
  where an integration by parts shows that
  \bas
	&&c_1 \xi^{n+\frac{(p-1)\alpha}{\beta}} \int_{\xi_1}^\xi \sigma^{-\frac{(p-1)\alpha}{\beta}-2} f'(\sigma) d\sigma
	= \Big(\frac{(p-1)\alpha}{\beta}+2\Big) c_1 \xi^{n+\frac{(p-1)\alpha}{\beta}}
	\int_{\xi_1}^\xi \sigma^{-\frac{(p+1)\alpha}{\beta}-3} f(\sigma) d\sigma \\
	& &\quad +~ c_1 \xi^{n-2} f(\xi)
	- c_1 \xi^{n+\frac{(p-1)\alpha}{\beta}} \xi_1^{-\frac{(p-1)\alpha}{\beta}-2} f(\xi_1) 
	\ge - c_1 \xi^{n+\frac{(p-1)\alpha}{\beta}} \xi_1^{-\frac{(p-1)\alpha}{\beta}-2} f(\xi_1) 
	\qquad \mbox{for all } \xi>\xi_1,
  \eas
  because $p\ge 1$ and $f$ is nonnegative.
  Using (\ref{56.4}), from (\ref{56.5}) we thus infer that
  \bea{56.6}
	g(\xi)
	&\ge& \Big\{ \xi_1^{-n} g(\xi_1) - c_1 \xi_1^{-2} f(\xi_1) \Big\}  \xi_1^{-\frac{(p-1)\alpha}{\beta}}
	\xi^{n+\frac{(p-1)\alpha}{\beta}} 
	\ge c_2 \xi^{n+\frac{(p-1)\alpha}{\beta}}
	\qquad \mbox{for all } \xi>\xi_1
  \eea
  with $c_2:= \frac{1}{2} \xi_1^{-n-(p-1)\alpha/\beta}$.\\
  In order to derive (\ref{56.1}) from this, we once again make use of the downward monotonicity of $f$ in estimating
  \bas
	g(\xi) = \int_0^\xi \sigma^{n-1} f^{1-p}(\sigma)d\sigma
	\le \int_0^\xi \sigma^{n-1} f^{1-p}(\xi) d\sigma
	= \frac{\xi^n}{n} f^{1-p}(\xi)
	\qquad \mbox{for all } \xi>0,
  \eas
  which combined with (\ref{56.6}) shows that
  \bas
	f(\xi)
	\le \Big\{ \frac{n}{\xi^n}  c_2 \xi^{n+\frac{(p-1)\alpha}{\beta}} \Big\}^{-\frac{1}{p-1}}
	= (nc_2)^{-\frac{1}{p-1}} \xi^{-\frac{\alpha}{\beta}}
	\qquad \mbox{for all } \xi>\xi_1.
  \eas
  Along with the boundedness of $f$ on $[0,\xi_1]$, this proves (\ref{56.1}).
\qed
\subsection{Implications on the large time behaviour for solutions of (\ref{0}) with algebraic spatial decay}
Let us underline the outcome of the above construction:
\begin{prop}\label{prop54}
  Let $p>1$ and $\alpha \in (0,\frac{1}{p})$, and set
$	\beta:=\frac{1-p\alpha}{2}$. 
  Then for any $A>0$,
  \bas
	\ua(x,t):=t^{-\alpha} \f \Big(t^{-\beta} |x| \Big),
	\qquad x\in\R^n, \ t>0,
  \eas
  defines a positive classical solution of the equation $u_t=u^p\Delta u$ in $\R^n \times (0,\infty)$.
\end{prop}
\proof
  Writing $u:=\ua$ and $f:=\f$ and abbreviating $\xi:=t^{-\beta}|x|$ for $x\in\R^n$ and $t>0$, we only need to compute
$	u_t(x,t)=-\alpha t^{-\alpha-1} f(\xi) - \beta t^{-\alpha-1} \xi f'(\xi)$
  and
  \bas
	u^p(x,t)\Delta u(x,t)
	= t^{-(p+1)\alpha-2\beta} f^p(\xi) \Big(f''(\xi)+\frac{n-1}{\xi} f'(\xi)\Big)
  \eas
  to see that by the choice of $\beta$ and (\ref{0f}) we have
  \bas
	u_t-u^p\Delta u
	= t^{-\alpha-1}  \bigg\{ - \alpha f'(\xi)- \beta\xi f'(\xi)
	- f^p(\xi) \Big(f''(\xi)+\frac{n-1}{\xi}f'(\xi)\Big)\bigg\} 
	= 0.
  \eas
\qed
\proofc of Theorem~\ref{theo1000}.\quad
  We only need to collect Proposition~\ref{prop54}, Lemma~\ref{lem53} and
Lemma~\ref{lem56}.
\qed
We can thereby complete the proof of Theorem~\ref{theo2000}.\abs
\proofc of Theorem~\ref{theo2000}~(ii).\quad
  Given $\gamma>0$, we let $\alpha:=\frac{\gamma}{p\gamma+2}$, so that clearly $\alpha\in (0,\frac{1}{p})$.
  Moreover, specifying $\beta:=\frac{1-p\alpha}{2}>0$, we then take $L>0$ as provided by
Lemma~\ref{lem53} and fix $A>0$ 
  suitably large fulfilling
  \be{55.3}
	LA\ge C_1.
  \ee
  Then with $\ua$ taken from Proposition~\ref{prop54}, the function $\ou$ defined by
  \be{55.4}
	\ou(x,t):=\ua(x,t+1)
	\equiv (t+1)^{-\alpha} \f \Big( (t+1)^{-\beta} |x|\Big),
	\qquad x\in\R^n, \ t\ge 0,
  \ee
  evidently satisfies $\ou_t=\ou^p \Delta \ou$ in $\R^n\times (0,\infty)$ and 
  \bas
	\ou(x,0)
	= \f (|x|)
	\ge LA  (1+|x|)^{-\frac{\alpha}{\beta}}
	\qquad \mbox{for all } x\in\R^n
  \eas
  by Lemma~\ref{lem53}. Since according to our definitions of $\beta$ and $\alpha$ we have
  $\frac{\alpha}{\beta}=\gamma$,
  in view of (\ref{55.3}) and (\ref{2000.1}) this shows that
$	\ou(x,0) \ge C_1  (1+|x|)^{-\gamma} \ge u_0(x)$
for all $x\in\R^n$.
  By means of a comparison argument, we thus conclude that for any $R>0$, the solution $u_R$ of (\ref{0R}) satisfies
  $u_R\le \ou$ in $B_R \times [0,\infty)$, and that hence 
  \be{55.5}
	u(x,t)\le \ou(x,t) 
	\qquad \mbox{for all $x\in\R^n$ and $t\ge 0$,}
  \ee
  which by (\ref{55.4}) and Lemma~\ref{lem51} in particular yields that
  \bas
	\|u(\cdot,t)\|_{L^\infty(\R^n)}
	\le \|\ou(\cdot,t)\|_{L^\infty(\R^n)}
	= A(t+1)^{-\alpha}
	\qquad \mbox{for all } t\ge 0.
  \eas
  Again using the definition of $\alpha$, we thereby immediately obtain (\ref{2000.2}) for $q=\infty$.
  Apart from that, once more writing $f:=\f$ we see that (\ref{55.5}) implies that whenever 
  $q\in (\frac{n}{\gamma},\infty)$, then
  \bea{55.6}
	\int_{\R^n} u^q(x,t) dx
	&\le& \int_{\R^n} \ou^q(x,t) dx 
	= (t+1)^{-q\alpha} \int_{\R^n} f^q \Big( (t+1)^{-\beta} |x| \Big) dx \nn\\
	&=& (t+1)^{-q\alpha}  (t+1)^{n\beta} \int_{\R^n} f^q(|y|) dy
	\qquad \mbox{for all } t\ge 0.
  \eea
  Here, applying Lemma~\ref{lem56} we find $C_2>0$ such that
$	f(\xi) \le C_2 \xi^{-\gamma} $ for all $\xi\ge 1$,
  so that our assumption $q>\frac{n}{\gamma}$ asserts that
  \bas
	\int_1^\infty \xi^{n-1} f^q(\xi)d\xi
	\le C_2^q \int_1^\infty \xi^{n-1-q\gamma} d\xi
	= \frac{C_2^q}{n-q\gamma}
  \eas
  is finite. 
  As furthermore, again by definition of $\beta$ and $\alpha$,
  \bas
	q\alpha-n\beta
	= q\alpha-\frac{n(1-p\alpha)}{2}
	= \frac{q\gamma}{p\gamma+2} - \frac{n}{2} \Big(1-\frac{p\gamma}{p\gamma+2}\Big)
	= \frac{q\gamma-n}{p\gamma+2},
  \eas
  from (\ref{55.6}) we thus infer the existence of $C_3>0$ such that
  \bas
	\int_{\R^n} u^q(x,t)dx
	\le C_3 (t+1)^{-\frac{q\gamma-n}{p\gamma+2}}
	\qquad \mbox{for all } t\ge 0,
  \eas
  which for these finite values of $q$ is equivalent to (\ref{2000.2}).
\qed

\noindent
{\bf Acknowledgements.} The first author was supported in part by the Slovak Research and
Development Agency under the contract No. APVV-14-0378 and by the VEGA grant
1/0319/15.

\end{document}